\documentclass[11pt]{amsart}
\usepackage{amsmath, amsthm, amscd, amssymb}
\usepackage[all]{xy}

\newtheorem{thm}{Theorem}[section]
\newtheorem*{thm*}{Theorem}
\newtheorem{cor}[thm]{Corollary}
\newtheorem*{cor*}{Corollary}
\newtheorem{lem}[thm]{Lemma}
\newtheorem{prop}[thm]{Proposition}

\newtheorem*{con*}{Conjecture}
\newtheorem*{prob*}{Problem}

\theoremstyle{definition}
\newtheorem{defn}[thm]{Definition}

\theoremstyle{remark}
\newtheorem{rem}[thm]{Remark}

\newcommand{\A}{\mathcal{A}}
\newcommand{\B}{\mathcal{B}}
\newcommand{\C}{\mathcal{C}}
\newcommand{\D}{\mathcal{D}}
\newcommand{\E}{\mathcal{E}}

\newcommand{\F}{\mathcal{F}}

\newcommand{\G}{\mathcal{G}}

\newcommand{\s}{\mathcal{S}}

\newcommand{\bbA}{\mathbb{A}}
\newcommand{\bbZ}{\mathbb{Z}}

\newcommand{\bbR}{\mathbb{R}}

\newcommand{\bbK}{\mathbb{K}}
\newcommand{\bbL}{\mathbb{L}}
\newcommand{\bbH}{\mathbb{H}}
\newcommand{\bbN}{\mathbb{N}}
\newcommand{\bbS}{\mathbb{S}}

\begin{document}

\title{A Controlled Approach to the Isomorphism Conjecture}


\begin{abstract}
We use a \textsl{hocolim} approach to the Isomorphism Conjecture in $K$-Theory to
analyze the case of groups of the form $G\rtimes\bbZ$ and $G_1*_{G}G_2$. As an important corollary we prove that the isomorphism conjecture in $K$-Theory holds for a finitely generated free group.
\end{abstract}

\author[Daniel Juan-Pineda]{Daniel Juan-Pineda$^*$}
\address{Instituto de Matem\'aticas,
UNAM Campus Morelia, Apartado Postal 61-3 (Xangari), Morelia,
 Mi\-cho\-a\-c\'an, MEXICO 58089}
\email{daniel@matmor.unam.mx}

\author[Stratos Prassidis]{Stratos Prassidis$^{**}$}
\address{Department of Mathematics
Canisius College,
Buffalo, NY 14208, U.S.A.}
\email{prasside@canisius.edu}

\protect\thanks{$^*$Partially supported by CONACyT and DGAPA-UNAM research grants.}
\protect\thanks{
$^{**}$Partially supported by a Canisius College summer research grant.}
\maketitle

\section{Introduction}

The Isomorphism Conjecture is, currently, one of the most important tools
in calculating algebraic invariants that appear in classification problems
in Topology. It should be considered as an {\it induction} technique that 
evaluates geometrically important obstruction groups of a space 
from the virtually cyclic subgroups of the
fundamental group of the space.
The algebraic invariants relevant to the rigidity problems in topology are the 
$A$ and $K$-invariants, $L^{-{\infty}}$-invariants and pseudoisotopy invariants. They
all can be characterized as elements of homotopy groups of the corresponding spectrum. 

More specifically, let ${\G}_{\Gamma}$ be the class of virtually cyclic subgroups of a 
group $\Gamma$,
${\E}{\G}_{\Gamma}$ the classifying space for the class ${\G}_{\Gamma}$, and 
$$p: E{\Gamma}{\times}_{\Gamma}{\E}{\G}_{\Gamma} \to {\E}{\G}_{\Gamma}/{\Gamma} =
{\B}{\G}_{\Gamma}$$
the projection map. Let $\bbS$ be any of the $\bbA$, $\bbK$, 
${\bbL}^{-{\infty}}$ or the pseudoisotopy 
spectra.

\begin{con*}[Isomorphism Conjecture (\cite{fajo:is})]
The assembly map 
$${\bbH}.({\B}{\G}_{\Gamma}, {\bbS}(p)) \to {\bbS}(B{\Gamma})$$
induced by the commutative diagram:
$$\begin{CD}
E{\Gamma}{\times}_{\Gamma}{\E}{\G}_{\Gamma} @>>> B{\Gamma} \\
@V{p}VV @VVV \\
{\B}{\G}_{\Gamma} @>>> *
\end{CD}$$
is a homotopy equivalence.
\end{con*}

The homology spectrum that appears in the Isomorphism Conjecture can be computed,
via a spectral sequence, from the ${\bbS}$-groups of the virtually cyclic subgroups of $\Gamma$.

Farrell--Jones have proved the isomorphism conjecture for a large class of geometrically
significant groups and for the pseudoisotopy spectrum (\cite{fajo:is}, \cite{fajo:ri-ge}). 
That implies the isomorphism conjecture for lower $K$-groups.

In this paper, we will study the Isomorphism Conjecture using homotopy colimits. The
homology spectrum in the statement of the conjecture can be described as a homotopy
colimit. In this context, the Isomorphism Conjecture reflects the extend of how much the
factors mentioned above commute with certain homotopy colimits. 

Homotopy colimits in the Isomorphism Conjecture have been used before 
(\cite{anmu:ge}, \cite{copr},
\cite{mupr:wa}, \cite{ta}). In all the references the homotopy colimit machinery
encodes, in a homological/homotopical way, controlled problems. 

The Isomorphism Conjecture can also be  considered as a statement about a ``forget control'' map.
In this context, we will use Segal's Pushdown Construction (\cite{hovo}) 
in homotopy colimits to capture the geometric idea of change of control. 
This approached is applied to two cases, with ${\bbS} = {\bbK}$, the $K$-theory spectrum:
\begin{itemize}
\item Groups that admit an epimorphism to $\bbZ$ and the kernel satisfying the IC.
\item Amalgamated free products where each factor satisfies the IC.
\end{itemize}

In the first case (Section \ref{sec-circle}), 
we will describe the part of $K$-theory that is controlled over the
circle and the forget control map in this case. More precisely, let ${\Gamma} = G 
{\rtimes}_{\alpha}{\bbZ}$ where $G$ satisfies the $\bbS$-IC.
The main result in this case (Theorem \ref{thm-circle}) states that the homology spectrum 
of $B{\G}_{\Gamma}{\times}S^1$ is homotopy equivalent to the ``controlled'' part, over $S^1$,
of the theory, in analogy with the results in \cite{hupr:co} and \cite{mupr:wa} in $K$-theory. 
Let ${\bbK}_R$ denotes the functor that maps a space $X$ to  $\bbK(R[\pi_1(X)])$.

\begin{thm*}[Semidirect Product with $\bbZ$]
Let ${\Gamma} = G{\rtimes}_{\alpha}{\bbZ}$ such that $G$ satisfies the ${\bbK}_R$-IC. 
Let
$$q_{\Gamma}: E{\Gamma}{\times}_{\Gamma}{\E}{\Gamma} \to S^1{\times}{\B}{\Gamma}.$$
Then 
\begin{enumerate}
\item there is an exact sequence:
$${\dots} \to {\pi}_i(BG) \xrightarrow{1 - {\alpha}_*} {\pi}_i(BG) \to
H_i(S^1{\times}{\B}{\Gamma}, {\bbK}_R(q_{\Gamma})) \to
{\pi}_{i-1}(BG) \xrightarrow{1 - {\alpha}_*} {\pi}_{i-1}(BG) \dots$$
\item If $RG$ is regular coherent, then the assembly map
$$H_i({\B}{\Gamma}, {\bbK}_R(p_{\Gamma})) \to K_i(R{\Gamma})$$
is an epimorphism for all $i\in \bbZ$.
\end{enumerate}
\end{thm*}

In the special case when $\Gamma$ is a
trivial extension, we derive information about the cokernel of the forget control map i.e. the
Nil-part of the theory.
That leads to a proof of the IC for the $K$-theory reduced Nil-spectrum (Section 
\ref{sec-trivial}), under the assumption that groups of the form $G_0{\times}{\bbZ}$,
with $G_0$ virtually infinite group, satisfy the ${\bbK}_R$-IC.

\begin{thm*}[Product with $\bbZ$]
Let $G$ satisfy the Bundle ${\bbK}_R$-IC and, for every virtually infinite 
cyclic subgroup $G_0 < G$, $G_0{\times}{\bbZ}$ satisfies the bundle ${\bbK}_R$-IC. 
Let ${\Gamma} = G{\times}{\bbZ}$. Then
\begin{enumerate}
\item ${\Gamma}$ satisfies the Bundle ${\bbK}_R$-IC.
\item $G$ satisfies the Bundle ${\bbN}il_R$-IC
\end{enumerate}
\end{thm*} 
In the second case (Section \ref{sec-interval}), we prove that the part of the
theory that is controlled over the interval satisfies a Mayer-Vietoris property in analogy to the
result in \cite{mupr:wa}.

\begin{thm*}[Amalgamated free products]
Let ${\Gamma} = G_1*_{G_0}G_2$ such that $G_i$, $i = 0, 1, 2$,  satisfy the ${\bbK}_R$-IC. 
Let
$$q_{\Gamma}: E{\Gamma}{\times}_{\Gamma}{\E}{\Gamma} \to I{\times}{\B}{\Gamma}.$$
Then 
\begin{enumerate}
\item there is an exact sequence:
$${\pi}_i(BG_0) \to {\pi}_i(BG_1){\oplus}{\pi}_i(BG_2) \to
H_i(I{\times}{\B}{\Gamma}, {\bbK}_R(q_{\Gamma})) \to
{\pi}_i(BG_0) \to {\pi}_{i-1}(BG_1){\oplus}{\pi}_{i-1}(BG_2)$$
\item If $RG_0$ is regular coherent, then the assembly map
$$H_i({\B}{\Gamma}, {\bbK}_R(p_{\Gamma})) \to K_i(R{\Gamma})$$
is an epimorphism for all $i\in\bbZ$.
\end{enumerate}
\end{thm*}

If in addition the {\it base groups} ($G$ in the semidirect product case and $G_0$ 
in the amalgamated free product case) are torsion free and satisfy the integral Novikov conjecture, then the  assembly map is a homotopy equivalence. Using these
ideas we prove (Corollary \ref{cor-free}):
 
\begin{cor*}
Let $F$ be a finitely generated free group and $R$ a regular coherent ring. Then 
\begin{enumerate}
\item $F$ satisfies the ${\bbK}_R$-IC.
\item $F{\rtimes}\bbZ$ satisfies the ${\bbK}_R$-IC
\end{enumerate}
\end{cor*}

The authors would like to thank Tom Farrell whose suggestions improved considerably the
original form of the paper. The second author would like to thank the Instituto
de Matem\'aticas, UNAM, Unidad Morelia for its hospitality during the preparation of this paper.

\section{Notation}\label{sec-not}

Let $\Gamma$ be a discrete group and ${\C}_{\Gamma}$ 
a class of subgroups of $\Gamma$ (i.e. a collection
of subgroups of $\Gamma$ closed under taking conjugates and subgroups). The {\it classifying
space} of the class ${\C}_{\Gamma}$, ${\E}{\C}_{\Gamma}$ is the $\Gamma$-complex 
whose isotropy groups are in ${\C}_{\Gamma}$
and its non-empty fixed point sets are contractible. A model for this space,
reminiscent of the bar construction (\cite{fajo:is}), is given as follows:
It is the realization of a semi-simplicial complex with 
$n$-simplex given by a sequence
$${\sigma} = 
{\gamma}_0{\Gamma}_0({\gamma}_1{\Gamma}_1, {\gamma}_2{\Gamma}_2, \dots , {\gamma}_n{\Gamma}_n)$$
with ${\Gamma}_i{\in}{\C}$ such that 
${\gamma}_i^{-1}{\Gamma}_{i-1}{\gamma}_i \subset {\Gamma}_i$ for
$i = 1, \dots , n$. The face operator is given by
$${\partial}_i{\sigma} = \left\{
\begin{array}{lc}
{\gamma}_0{\gamma}_1{\Gamma}_1({\gamma}_2{\Gamma}_1, \dots , {\gamma}_n{\Gamma}_n), & i = 0, \\
{\gamma}_0{\Gamma}_0({\gamma}_1{\Gamma}_1, \dots , 
{\gamma}_{i-1}{\Gamma}_{i-1}, {\gamma}_i{\gamma}_{i+1}{\Gamma}_{i+1}, 
{\gamma}_{i+1}{\Gamma}_{i+1}, \dots ,
{\gamma}_n{\Gamma}_n), & 0 < i < n ,\\
{\gamma}_0{\Gamma}_0({\gamma}_1{\Gamma}_1, \dots , {\gamma}_{n-1}{\Gamma}_{n-1}), & i = n.
\end{array}
\right.
$$
The group $\Gamma$ acts on $E{\C}_{\Gamma}$ by
$${\gamma}{\sigma} = ({\gamma}{\gamma}_0){\Gamma}_0({\gamma}_1{\Gamma}_1, 
{\gamma}_2{\Gamma}_2, \dots , {\gamma}_n{\Gamma}_n), \; 
\text{for} \; {\gamma}{\in}{\Gamma}.$$
We write ${\B}{\C}_{\Gamma}$ for the orbit space ${\E}{\C}_{\Gamma}/{\Gamma}$. 
To ensure that ${\B}{\C}_{\Gamma}$
is a simplicial complex, we subdivide ${\E}{\C}_{\Gamma}$ twice, i.e.
${\B}{\C}_{\Gamma} = ({\E}{\C}_{\Gamma})''/{\Gamma}$.

The construction of the classifying space is functorial with respect to group homomorphisms.
Let ${\C}_{\Gamma}$ be a class of subgroups of $\Gamma$.
For if ${\rho}: {\Gamma} \to G$ be a group homomorphism then $\rho$ induces a 
$\rho$-equivariant map
$$\bar{\rho}: {\E}{\C}_{\Gamma} \to {\E}{\C}_G, \;
{\gamma}_0{\Gamma}_0({\gamma}_1{\Gamma}_1, {\gamma}_2{\Gamma}_2, \dots , {\gamma}_n{\Gamma}_n)
\mapsto
{\rho}({\gamma}_0){\rho}({\Gamma}_0)
({\rho}({\gamma}_1){\rho}({\Gamma}_1), 
{\rho}({\gamma}_2){\rho}({\Gamma}_2), \dots , {\rho}({\gamma}_n){\rho}({\Gamma}_n))
$$
where ${\C}_G$ is a class of subgroups of $G$ that contains the images, under $\rho$, of the
elements of ${\C}_{\Gamma}$.
The map $\bar{\rho}$ induces a map ${\rho}'$ to the quotient spaces.

For each simplicial complex $K$, we write $\text{cat}(K)$ for the category of simplices of
$K$, viewed as a partially ordered set. 
Thus objects are the simplices of $K$ and there is a single morphism
from $\sigma$ to $\tau$ whenever ${\sigma} \le {\tau}$. 

\begin{defn} [\cite{anmu:ge}]
Let $p: E \to B$ be a map with $B = |K|$, the geometric realization of a simplicial complex. 
The map is said to have a {\it homotopy colimit structure} if there is a functor 
$$F: \text{cat}(K)^{op} \to {\bf Top}$$
such that:
\begin{itemize}
\item $E = \text{hocolim}_{\text{cat}(K)^{op}}(F)$.
\item $p = \text{hocolim}_{\text{cat}(K)^{op}}({\nu})$, where $\nu$ is the natural transformation
from $F$ to the constant point functor.
\end{itemize}
\end{defn}

\noindent
\underline{Notation}. Let $p: E \to B$ be a map, $B = |K|$ and $\text{cat}(K)^{op}$ the category
of the simplicial complex $K$. We define the {\it barycentre functor}
$$\text{bar}(p): \text{cat}(K)^{op} \to {\bf Top}, \;
{\sigma} \mapsto p^{-1}(\hat{\sigma})$$
where $\hat{\sigma}$ is the barycentre of $\sigma$.

\vspace{12pt}

\begin{rem}\label{rem-hocolim} 
We give basic examples of maps that admit a homotopy colimit structure. The proofs follow
from direct calculations (also \cite{ta}).
\begin{enumerate}
\item Let $p: E \to B$, $B = |K|$, be a map that admits a homotopy colimit structure
relative to the functor $\text{bar}(p)$. Let $p': E' \to E$ be a fiber bundle. Then
the composition
$$q: E' \xrightarrow{p'} E \xrightarrow{p} B$$
admits a homotopy colimit structure relative to the functor $\text{bar}(q)$.
\item Let $\Gamma$ be a discrete group, $E{\Gamma}$ a free contractible $\Gamma$-complex and
${\E}{\C}_{\Gamma}$ the classifying complex for the family of subgroups of
$\Gamma$. Let
$$p_{\Gamma}: E{\Gamma}{\times}_{\Gamma}{\E}{\C}_{\Gamma} \to {\E}{\C}_{\Gamma}/{\Gamma} = 
{\B}{\C}_{\Gamma}$$
be the projection map to the second coordinate. Then $p_{\Gamma}$ 
has a homotopy colimit structure with respect to the functor
$\text{bar}(p_{\Gamma})$. Notice that, in this case,
$p_{\Gamma}^{-1}(\hat{\sigma})$ is a space of type 
$B{\Gamma}_{\sigma}$, where ${\Gamma}_{\sigma}$ is the isotropy group of $\sigma$, an
element in the class ${\C}_{\Gamma}$.
\item Let ${\rho}: {\Gamma} \to G$ be a group epimorphism. Then the map
$$q: E{\Gamma}{\times}_{\Gamma}{\E}{\C}_{\Gamma} \;\xrightarrow{p}\;
{\B}{\C}_{\Gamma} \;\xrightarrow{{\rho}'}\; {\B}{\C}_G$$
has a homotopy colimit structure with respect to the functor
$\text{bar}(q)$, where ${\C}_G = {\rho}({\C}_{\Gamma})$, the class of subgroups of $G$ 
consisting of the images of elements of ${\C}_{\Gamma}$.
\end{enumerate}
\end{rem}

Let $F: {\C} \to {\D}$ and $X: {\C} \to {\bf Top}$ be two functors. Then {\it
Segal's Pushdown Construction} (see for example \cite{hovo})
defines a functor $F_*X: {\D} \to {\bf Top}$ such that
$$\text{hocolim}_{\C}X \simeq \text{hocolim}_{\D}F_*X.$$
We will explicitely describe the construction to the case of Part(3) in Remark
\ref{rem-hocolim}. In this case, we start with a map:
$$q: E{\Gamma}{\times}_{\Gamma}{\E}{\C}_{\Gamma} \;\xrightarrow{p}\;
{\B}{\C}_{\Gamma} \;\xrightarrow{{\rho}'}\; {\B}{\C}_G$$
The map ${\rho}'$ induces a functor
$$P: \text{cat}({\B}{\C}_{\Gamma})^{op} \to \text{cat}({\B}{\C}_G)^{op}$$
We will describe the functor
$$P_*\text{bar}(p): \text{cat}({\B}{\C}_G)^{op} \to {\bf Top}.$$
For each simplex ${\sigma}$ of ${\B}{\C}_G$, let $P{\downarrow}{\sigma}$ be the over
category. In this case, the objects of $P{\downarrow}{\sigma}$ are simplices $\tau$
of ${\B}{\C}_{\Gamma}$ such that ${\rho}'({\tau})$ contains $\sigma$ as a face. Set
$$p_{\sigma} = p|: p^{-1}(|P{\downarrow}{\sigma}|) \to |P{\downarrow}{\sigma}|.$$
Then $P_*\text{bar}(p)({\sigma}) = 
\text{hocolim}_{P{\downarrow}{\sigma}}\text{bar}(p_{\sigma})$.
Summarizing:

\begin{prop}\label{prop-hoco}
There is a homotopy equivalence:
$$\text{hocolim}_{\text{cat}({\B}{\C}_{\Gamma})^{op}}(\text{bar}(p)) \simeq
\text{hocolim}_{\text{cat}({\B}{\C}_G)^{op}}(P_*\text{bar}(p)).$$
\end{prop}

The class of subgroups of $\Gamma$ of interest in the Isomorphism Conjecture is the class of
{\it virtually cyclic subgroups}, denoted ${\G}_{\Gamma}$. They naturally split into 
two categories:
\begin{itemize}
\item Finite subgroups of $\Gamma$.
\item Virtually infinite cyclic subgroups of $\Gamma$ i.e. subgroups which contain 
an infinite cyclic subgroup of finite index.
\end{itemize}

The subgroups of the second type are two-ended subgroups of $\Gamma$ (\cite{dd}) and they split 
into two types:
\begin{itemize}
\item Groups $H$ that admit an epimorphism to $\bbZ$ with finite kernel i.e. 
  $$H \cong K{\rtimes}{\bbZ}.$$
\item Groups $H$ that admit an epimorphism to the infinite dihedral subgroup $D_{\infty}$ with
finite kernel, i.e.
$$H \cong A*_BC, \quad [B : 1] < {\infty}, \; [A : B] = [C : B] = 2.$$
\end{itemize}

\section{A Homotopy Approach to the Isomorphism Conjecture}

We will present a reformulation of the IC using homotopy colimits. This approach depends
heavily on the interpretation of the homology spectrum as a homotopy colimit (\cite{anmu:ge}, 
\cite{ta}).

Let $\Gamma$ be a discrete group and $\bbS$ a homotopy invariant functor from spaces to spectra.
As before, let ${\G}_{\Gamma}$ be the class of virtually cyclic (finite or infinite)
subgroups of $\Gamma$. Let ${\E}{\G}_{\Gamma}$ be the classifying $\Gamma$-complex
for the class ${\G}_{\Gamma}$ and
$$p_{\Gamma}: E{\Gamma}{\times}_{\Gamma}{\E}{\G}_{\Gamma} \to {\E}{\G}_{\Gamma}/{\Gamma} = 
{\B}{\G}_{\Gamma}$$
be the projection map. Let $\text{cat}({\B}{\G}_{\Gamma})^{op}$ be the
category corresponding to the partially ordered set of simplices of ${\B}{\G}_{\Gamma}$.
Let $\text{bar}(p_{\Gamma})$ be the barycentre functor. 
Let $r: Y \to B{\Gamma}$ be
a bundle. Form the pull-back:
$$\begin{CD}
\bar{Y} @>>> Y \\
@V{\rho}VV @VV{r}V \\
E{\Gamma}{\times}_{\Gamma}{\E}{\G}_{\Gamma} @>>> B{\Gamma} \\
@V{p_{\Gamma}}VV @VVV \\
{\B}{\G}_{\Gamma} @>>> *
\end{CD}$$
Notice that
$$\text{hocolim}_{\text{cat}({\B}{\G}_{\Gamma})^{op}}\text{bar}(p_{\Gamma}{\circ}{\rho}) \cong
\bar{Y} \simeq Y$$
Let $*$ denote the category with a single object and a single morphism. Let
$$F_{\Gamma}: * \to {\bf Top}, \quad F_{\Gamma}(*) = Y$$

\vspace{18pt}
\noindent
\underline{\it The Bundle Isomorphism Conjecture}. With the notation above, the functor from 
$\text{cat}({\B}{\G}_{\Gamma})^{op}$ to $*$, induces a homotopy equivalence of 
spectra:
$$\text{hocolim}_{\text{cat}({\B}{\G}_{\Gamma})^{op}}
{\bbS}{\circ}\text{bar}(p_{\Gamma}{\circ}{\rho}) \to
\text{hocolim}_{\text{cat}(*)}{\bbS}{\circ}\text{bar}(F_{\Gamma}) \cong
{\bbS}(Y).$$

\vspace{18pt}
\begin{rem} \hfill

\begin{enumerate}
\item In the interesting cases, the functor $\bbS$ factors as
$${\bbS}: {\bf Top} \xrightarrow{\s} {\C} \xrightarrow{B} {\bbS}pectra$$
where $\C$ is the category of small categories and $B$ is the classifying spectrum functor. 
Using the notation in \cite{th:fi}, Theorem 3.19,
the left hand side can also be re-written:
$$\begin{array}{lrl}
\text{hocolim}_{\text{cat}({\B}{\G}_{\Gamma})^{op}}
B({\s}{\circ}\text{bar}(p_{\Gamma}{\circ}{\rho})) & \simeq & 
B\left(\displaystyle{\int_{{\text{cat}({\B}{\G}_{\Gamma})^{op}}}
{\s}{\circ}\text{bar}(p_{\Gamma}{\circ}{\rho})}
\right) \\ [2ex]
& &
\end{array}$$
Thus the isomorphism conjecture can be reformulated as follows:

\vspace{18pt}
\noindent
\underline{\it Categorical Reformulation of the Bundle IC}. The natural functor 
$\text{cat}({\B}{\G}_{\Gamma})^{op} \to *$ induces a homotopy equivalence of categories:
$$\int_{{\text{cat}({\B}{\G}_{\Gamma})^{op}}}{\s}{\circ}\text{bar}(p_{\Gamma}{\circ}{\rho}) \to 
\int_*{\s}{\circ}F_{\Gamma} \;\cong\; {\s}(Y)$$
\item Since $\bbS$ is homotopy invariant and 
$Y \simeq \text{hocolim}_{\text{cat}({\B}{\G}_{\Gamma})^{op}}\text{bar}(p_{\Gamma}{\circ}{\rho})$
we see that the Bundle IC states that $\s$ commutes,
up to homotopy, with the homotopy colimit above, i.e.:

\vspace{18pt}\noindent
\underline{\it IC Using Homotopy Colimits}.
With the above notation,
$$\int_{{\text{cat}({\B}{\G}_{\Gamma})^{op}}}{\bbS}
{\circ}\text{bar}(p_{\Gamma}{\circ}{\rho}) \simeq 
{\bbS}(\text{hocolim}_{\text{cat}({\B}{\G}_{\Gamma})^{op}}\text{bar}(p_{\Gamma}{\circ}{\rho})).$$
\item 
The above formulation is related to the classical statement of the conjecture in \cite{fajo:is}
as follows:
\begin{enumerate}
\item If $r$ is the identity map, then the Bundle IC is the Isomorphism Conjecture in 
\cite{fajo:is}.
\item In \cite{fajo:is}, the Fibered Isomorphism Conjecture is stated under the assumption that
$r$ is a fibration. Our assumption that $r$ is a bundle is needed to preserve the homotopy
colimit structure of the composition $p_{\Gamma}{\circ}{\rho}$. 
\end{enumerate}
\end{enumerate}
\end{rem}

As an illustration of the categorical approach to the Isomorphism Conjecture, we
study the Isomorphism Conjecture under certain group extensions. 

Let $\bbS$ be
a homotopy invariant functor from spaces to spectra and ${\rho}: {\Gamma} \to G$ a
group epimorphism. Then the map $\rho$ induces a commutative diagram (Section \ref{sec-not}):
$$\begin{CD}
E{\Gamma}{\times}_{\Gamma}{\E}{\G}_{\Gamma} @>>> EG{\times}_G{\E}{\G}_G\\
@V{p_{\Gamma}}VV @VV{p_G}V \\
{\B}{\G}_{\Gamma} @>{{\rho}'}>> {\B}{\G}_G 
\end{CD}$$
The top horizontal map is induced by the map between classifying spaces and the map
$$\bar{\rho}: {\E}{\G}_{\Gamma} \to {\E}{\G}_G.$$
The map $p_{\Gamma}$ induces a functor:
$$P_{\Gamma}: \text{cat}({\B}{\G}_{\Gamma})^{op} \to \text{cat}({\B}{\G}_G)^{op}.$$

Let $[\sigma]_G$ represent a simplex of ${\B}{\G}_G$.
The space ${\E}{\G}_G$ with
the $G_{\sigma}$ action is a space of type ${\E}{\G}_{G_{\sigma}}$ and
the quotient ${\E}{\G}_G/G_{\sigma}$ is a space of type ${\B}{\G}_{G_{\sigma}}$.
Since $G_{\sigma}$ is virtually cyclic, ${\E}{\G}_G$ is $G_{\sigma}$-contractible and 
${\B}{\G}_{G_{\sigma}}$ is contractible.
We write 
$$p_{G_{\sigma}}: EG{\times}_{G_{\sigma}}{\E}{\G}_G\to {\B}{\G}_{G_{\sigma}}$$
for the projection map. 

Let ${\Delta}_{\sigma} = {\rho}^{-1}(G_{\sigma})$, a subgroup of $\Gamma$. 
Then ${\E}{\G}_{\Gamma}$,
with the ${\Delta}_{\sigma}$ action, is a space of type ${\E}{\G}_{{\Delta}_{\sigma}}$. Consider
the projection map:
$$p_{{\Delta}_{\sigma}}: E{\Gamma}{\times}_{{\Delta}_{\sigma}}{\E}{\G}_{\Gamma} \to 
{\E}{\G}_{\Gamma}/{\Delta}_{\sigma} = {\B}{\G}_{{\Delta}_{\sigma}}$$
Then the functor $\text{bar}(p_{{\Delta}_{\sigma}})$ 
induces the homotopy colimit structure on $p_{{\Delta}_{\sigma}}$. The map $\rho$ induces a
commutative diagram:
$$\begin{CD}
E{\Gamma}{\times}_{{\Delta}_{\sigma}}{\E}{\G}_{\Gamma} @>>> EG{\times}_{G_{\sigma}}{\E}{\G}_G \\
@V{p_{{\Delta}_{\sigma}}}VV @VV{p_{G_{\sigma}}}V \\
{\B}{\G}_{{\Delta}_{\sigma}} @>{{\rho}_{\sigma}}>> {\B}{\G}_{G_{\sigma}}
\end{CD}$$ 
With this set up, the map ${\rho}_{\sigma}$ maps simplices to simplices and induces a functor
$$P_{\sigma}: \text{cat}({\B}{\G}_{{\Delta}_{\sigma}})^{op} \to 
\text{cat}({\B}{\G}_{G_{\sigma}})^{op}$$

\begin{lem}\label{lem-over}
For each simplex $[{\sigma}]$ of ${\B}{\G}_G$ there are homeomorphisms:
$${\phi}_{\sigma}: |P_{\sigma}{\downarrow}[{\sigma}]_G| \xrightarrow{\cong} 
|P_{\Gamma}{\downarrow}[{\sigma}]_G|,
\quad
{\Phi}_{\sigma}: p_{{\Delta}_{\sigma}}^{-1}|P_{\sigma}{\downarrow}[{\sigma}]_G| 
\xrightarrow{\cong}
p_{\Gamma}^{-1}|P_{\Gamma}{\downarrow}[{\sigma}]_G|$$
natural in $\sigma$.
\end{lem}

\begin{proof}
We will show that the first map is a homeomorphism. 
The inclusion map ${\Delta}_{\sigma} \to \Gamma$ induces  a map
$${\phi}: {\B}{\G}_{{\Delta}_{\sigma}} = {\E}{\G}_{\Gamma}/{\Delta}_{\sigma} \to
{\E}{\G}_{\Gamma}/{\Gamma} = {\B}{\G}_{\Gamma}.$$
The restriction of $\phi$ induces a map 
$${\phi}_{\sigma}: |P_{\sigma}{\downarrow}[{\sigma}]_G| \to
|P_{\Gamma}{\downarrow}[{\sigma}]_G|$$
We will define the inverse of ${\phi}_{\sigma}$. For this, let $[{\tau}]_{\Gamma}$ represent
a simplex in $|P_{\Gamma}{\downarrow}[{\sigma}]_G|$. Thus there is $g \in G$ such
that $g\bar{{\rho}}({\tau}) > {\sigma}$. 
Let ${\gamma} \in {\Gamma}$ be such that ${\rho}({\gamma}) 
= g$. Define
$${\psi}_{\sigma}: |P_{\Gamma}{\downarrow}[{\sigma}]_G| \to |P_{\sigma}{\downarrow}[{\sigma}]_G|,
\quad {\psi}_{\sigma}([{\tau}_{\Gamma}]) = [{\gamma}{\tau}]_{{\Delta}_{\sigma}}$$
\begin{itemize}
\item ${\psi}_{\sigma}$ is well defined:
\begin{itemize}
\item ${\psi}_{\sigma}([{\tau}_{\Gamma}]) \in |P_{\sigma}{\downarrow}[{\sigma}]_G|$: That
follows because 
$$\bar{\rho}({\gamma}{\tau}) = {\rho}({\gamma})\bar{\rho}({\tau}) =
g\bar{\rho}({\tau}) > {\sigma}.$$
\item The definition of ${\psi}_{\sigma}$ does not depend on the choice of $\gamma$: Let
${\gamma}'\in \Gamma$ be such that ${\rho}({\gamma}') = g$. Then
$${\gamma}^{-1}{\gamma}' \in \text{Ker}({\rho}) \subset {\Delta}_{\sigma}.$$
Thus $[{\gamma}{\tau}]_{{\Delta}_{\sigma}} = [{\gamma}'{\tau}]_{{\Delta}_{\sigma}}$.
\item The definition of ${\psi}_{\sigma}$ does not depend on the choice 
of the representative of the orbit of $\tau$: Let $g'\in G$ such that $g'\bar{\rho}({\tau}) > 
{\sigma}$. Then $\bar{\rho}({\tau})$ contains both $g^{-1}{\sigma}$ and $(g')^{-1}{\sigma}$.
Since ${\B}{\G}_G$ is a simplicial complex, that implies 
$$g^{-1}{\sigma} = (g')^{-1}{\sigma} \;\Rightarrow\; g'g^{-1} \in G_{\sigma}.$$
Since ${\Delta}_{\sigma} = {\rho}^{-1}(G_{\sigma})$, ${\psi}_{\sigma}$ is well defined.
\end{itemize}
\item ${\psi}_{\sigma}$ is the inverse of ${\phi}_{\sigma}$: That follows directly from the
definition of the functions.
\end{itemize}

For the construction of ${\Phi}_{\sigma}$ 
we start by defining a subcomplex of ${\E}{\G}_{\Gamma}$: 
$$Q_{\sigma} = \{{\tau}\in {\E}{\G}_{\Gamma}, \; \bar{\rho}({\tau}) > {\sigma}\}$$
Then
$$p_{{\Delta}_{\sigma}}^{-1}(|P_{\sigma}{\downarrow}{\sigma}|) = 
E{\Gamma}{\times}_{{\Delta}_{\sigma}}{\Delta}_{\sigma}Q_{{\sigma}}, \quad 
p_{\Gamma}^{-1}(|P_{\Gamma}{\downarrow}{\sigma}|) = 
E{\Gamma}{\times}_{\Gamma}{\Gamma}Q_{\sigma}.$$
Since each simplex in $Q_{\sigma}$ has isotropy group a subgroup of ${\Delta}_{\sigma}$,
$$p_{\Gamma}^{-1}(|P_{\Gamma}{\downarrow}{\sigma}|) = 
E{\Gamma}{\times}_{\Gamma}{\Gamma}Q_{\sigma} \cong 
E{\Gamma}{\times}_{{\Delta}_{\sigma}}{\Delta}_{\sigma}Q_{\sigma}.$$
The naturality of the homeomorphisms is immediate from their construction.
\end{proof}

As in Proposition \ref{prop-hoco}, there is a homotopy equivalence
$$\text{hocolim}_{\text{cat}({\B}{\G}_{\Gamma})^{op}}({\bbS}{\circ}\text{bar}(p_{\Gamma})) \simeq
\text{hocolim}_{\text{cat}({\B}{\G}_G)^{op}}(P_{{\Gamma}*}({\bbS}{\circ}\text{bar}(p_{\Gamma}))).$$
Using Lemma \ref{lem-over}, we can reformulate the description of the
functor $P_{{\Gamma}*}({\bbS}{\circ}\text{bar}(p_{\Gamma}))$. For this, we define
$$h_{\Gamma}: \text{cat}({\B}{\G}_G)^{op} \to {\bbS}pectra, \;
{\sigma} \mapsto \text{hocolim}_{\text{cat}({\B}{\G}_{{\Delta}_{\sigma}})^{op}}({\bbS}{\circ}
\text{bar}(p_{{\Delta}_{\sigma}})).$$

\begin{prop}\label{prop-reduction}
With the above notation, there is a homotopy equivalence
$$\text{hocolim}_{\text{cat}({\B}{\G}_{\Gamma})^{op}}({\bbS}{\circ}\text{bar}(p_{\Gamma})) \simeq
\text{hocolim}_{\text{cat}({\B}{\G}_G)^{op}}(h_{\Gamma}).$$
\end{prop}

\begin{proof}
We put everything together:
\begin{itemize}
\item Using the commutative diagram for $G_{\sigma}$ and Segal's Pushdown Construction,
$$\text{hocolim}_{\text{cat}({\B}{\G}_{{\Delta}_{\sigma}})^{op}}
({\bbS}{\circ}\text{bar}(p_{{\Delta}_{\sigma}}))  \simeq 
\text{hocolim}_{\text{cat}({\B}{\G}_{G_{\sigma}})^{op}}(P_{{\sigma}*}({\bbS}{\circ}
\text{bar}(p_{{\Delta}_{\sigma}})))$$
\item Since ${\B}{\G}_{G_{\sigma}}$ contracts to $[\sigma]_G$, which is also contractible,
$$\text{hocolim}_{\text{cat}({\B}{\G}_{{\Delta}_{\sigma}})^{op}}
({\bbS}{\circ}\text{bar}(p_{{\Delta}_{\sigma}}))  \simeq 
P_{{\sigma}*}({\bbS}{\circ}\text{bar}(p_{{\Delta}_{\sigma}}))([{\sigma}]_G) \simeq \text{
hocolim}_{P_{\sigma}{\downarrow}[{\sigma}]_G}({\bbS}{\circ}\text{bar}(p_{{\Delta}_{\sigma}}|)).$$
\item By Lemma \ref{lem-over}, the diagram
$$\begin{CD}
p_{{\Delta}_{\sigma}}^{-1}(|P_{\sigma}{\downarrow}[{\sigma}]_G|) @>>> 
p_{\Gamma}^{-1}(|P_{\Gamma}{\downarrow}[{\sigma}]_G|)\\
@V{p_{{\Delta}_{\sigma}}}VV @VV{p_{\Gamma}}V \\
|P_{\sigma}{\downarrow}[{\sigma}]_G| @>>> |P_{\Gamma}{\downarrow}[{\sigma}]_G|
\end{CD}$$
commutes and the horizontal maps are homeomorphisms. Thus, by naturality,
$$ \text{
hocolim}_{P_{\sigma}{\downarrow}[{\sigma}]_G}({\bbS}{\circ}\text{bar}(p_{{\Delta}_{\sigma}})|)
\cong
\text{
hocolim}_{P_{\Gamma}{\downarrow}[{\sigma}]_G}({\bbS}{\circ}\text{bar}(p_{\Gamma})|).$$
\item Thus
$$h_{\Gamma}({\sigma}) = 
\text{hocolim}_{\text{cat}({\B}{\G}_{{\Delta}_{\sigma}})^{op}}({\bbS}{\circ}
\text{bar}(p_{{\Delta}_{\sigma}})) \simeq
\text{hocolim}_{P_{\Gamma}{\downarrow}[{\sigma}]_G}({\bbS}{\circ}\text{bar}
(p_{\Gamma})|)$$
and the equivalence is natural in $\sigma$.
\item Since $P_{{\Gamma}*}\text{bar}(p_{\Gamma})({\sigma}) =
 \text{hocolim}_{P_{\Gamma}{\downarrow}[{\sigma}]_G}({\bbS}{\circ}\text{bar}
(p_{\Gamma})|)$,
$$\text{hocolim}_{\text{cat}({\B}{\G}_{\Gamma})^{op}}({\bbS}{\circ}\text{bar}(p_{\Gamma})) \simeq 
\text{hocolim}_{\text{cat}({\B}{\G}_G)^{op}}(P_{{\Gamma}*}\text{bar}(p_{\Gamma})) \simeq
\text{hocolim}_{\text{cat}({\B}{\G}_G)^{op}}(h_{\Gamma})$$
\end{itemize}
That completes the proof of the proposition.
\end{proof}

The next Lemma is an application of Proposition \ref{prop-hoco} to the barycentre functors.

\begin{lem}\label{lem-point}
With the above notation,  there is a homotopy equivalence
$$B{\Gamma} \simeq \text{hocolim}_{\text{cat}({\B}{\G}_G)^{op}}({\beta})$$
where ${\beta}({\sigma}) = \text{hocolim}_{P_{\Gamma}{\downarrow}[{\sigma}]_G}
\text{bar}(p_{\Gamma}|)$ is
a space of type $B{\Delta}_{\sigma}$.
\end{lem}

\begin{proof}
Notice that a model for $B{\Gamma}$ is given by 
$\text{hocolim}_{\text{cat}({\B}{\G}_{\Gamma})^{op}}(\text{bar}(p_{\Gamma}))$. 
The result follows from 
applying Proposition \ref{prop-hoco}.
\end{proof}

\begin{thm}\label{thm-projection}
Let $\bbS$ be a homotopy invariant functor from spaces to spectra and ${\rho}: {\Gamma} \to G$
a group epimorphism. Assume that
\begin{enumerate}
\item The Bundle $\bbS$-IC is true for $G$.
\item For each $F{\in}{\G}_G$, the Bundle $\bbS$-IC is true
for ${\rho}^{-1}(F)$.
\end{enumerate}
Then the Bundle $\bbS$-IC holds for $\Gamma$.
\end{thm}

\begin{proof}
We will give the proof when the bundle over $B{\Gamma}$ is the identity. The general case
follows similarly.
By Proposition \ref{prop-reduction},
$$\text{hocolim}_{\text{cat}({\B}{\G}_{\Gamma})^{op}}({\bbS}{\circ}\text{bar}(p_{\Gamma})) \simeq
\text{hocolim}_{\text{cat}({\B}{\G}_G)^{op}}(h_{\Gamma})$$
where 
$$h_{\Gamma}([{\sigma}]_G) = 
\text{hocolim}_{\text{cat}({\B}{\G}_{{\Delta}_{\sigma}})^{op}}({\bbS}{\circ}
\text{bar}(p_{{\Delta}_{\sigma}}))$$
with ${\Delta}_{\sigma} = {\rho}^{-1}(G_{\sigma})$. By Assumption (2), the functor $h_{\Gamma}$ is
naturally homotopy equivalent to the functor $h'$, given by
$$h'([{\sigma}]_G) = {\bbS}(B{\Delta}_{\sigma}) \simeq {\bbS}({\beta}([{\sigma}]_G))$$
The group $\Gamma$ acts on $EG$ via $\rho$. The action is not free because elements of the 
kernel fix $EG$. Consider the bundle map:
$$B{\Gamma} = E{\Gamma}{\times}_{\Gamma}EG \to EG/G = BG$$
induced by projection to the second coordinate. Form the pull-back diagram:
$$\begin{CD}
Y @>>> B{\Gamma} \\
@V{r}VV @VVV \\
EG{\times}_G{\E}{\G}_G @>>> BG \\
@V{p_G}VV @VVV \\
{\B}{\G}_G @>>> *
\end{CD}$$
For each $[{\sigma}]_G$  a simplex in ${\B}{\G}_G$, the inverse image of its barycentre
under $p_G$ is a space of type $BG_{\sigma}$.  Thus the inverse image
$$(p_G{\circ}r)^{-1}([\hat{\sigma}]_G) \simeq B{\Delta}_{\sigma}.$$
Since the Bundle $\bbS$-IC holds for $G$
$$\begin{array}{lll}
{\bbS}(B{\Gamma}) &\simeq &
\text{hocolim}_{\text{cat}({\B}{\G}_G)^{op}}({\bbS}{\circ}\text{bar}(p_G{\circ}r)) \\
& \simeq & 
\text{hocolim}_{\text{cat}({\B}{\G}_G)^{op}}({\bbS}(B{\Delta}_{\sigma})) \\
& \simeq &
\text{hocolim}_{\text{cat}({\B}{\G}_G)^{op}}(h') \\
&\simeq & \text{hocolim}_{\text{cat}({\B}{\G}_G)^{op}}(h_{\Gamma}) \\
&\simeq & \text{hocolim}_{\text{cat}({\B}{\G}_{\Gamma})^{op}}({\bbS}{\circ}\text{bar}(p_{\Gamma}))
\end{array}
$$
proving the $\bbS$-IC for $\Gamma$.
\end{proof}

\begin{cor}\label{cor-bun}
The proof implies that
if the Bundle $\bbS$-IC holds for $G$ and the $\bbS$-IC holds for the
inverse images of virtually cyclic subgroups of $G$, then the $\bbS$-IC holds for $\Gamma$.
\end{cor}

\section{Spaces over the Circle}\label{sec-circle}

Let ${\Gamma} = G {\rtimes}_{\alpha}{\bbZ}$. Here $\alpha$ is the automorphism of
$G$ induced by the action of the generator of ${\bbZ}$. The automorphism $\alpha$
is well-defined up to inner automorphisms.
By choosing a suitable right $G$-space for $EG$, there is an $\alpha$-equivariant 
homeomorphism
$${\psi}: EG \to EG$$
i.e. ${\psi}(xg) = {\psi}(x){\alpha}(g)$. By taking quotients, we see that there is
a homeomorphism
$${\phi}: BG \to BG$$
that induces the map $\alpha$ in the fundamental group, again up to inner automorphisms.
Choose as a model for $B{\Gamma}$ the mapping torus of $\phi$:
$$B{\Gamma} = BG{\times}[0, 1]/{\sim}, \; ({\phi}(x), 0) \sim (x, 1).$$
Then a model for $E{\Gamma}$ is the infinite mapping telescope of ${\psi}$:
$$E{\Gamma} = EG{\times}[0, 1]{\times}{\bbZ}/{\sim}, \quad
(x, 1, n) \sim ({\psi}(x), 0, n + 1).$$
The right action of $\Gamma$ on $E{\Gamma}$ is given by:
$$(x, t, n)(g, m) = (x{\alpha}^n(g), t, n + m).$$
The map 
$${\Phi}: BG{\times}[0, 1] \to B{\Gamma}, \; {\Phi}(x, t) = [x, t]$$
has the property that:
$${\Phi}(x, 0) = [x, 0], \quad \Phi(x, 1) = [x, 1] = [{\phi}(x), 0]$$ 
i.e., it defines a homotopy between the identity map and the map $\phi$ on $BG$ inside $B{\Gamma}$.
Also, the natural projection map to the second coordinate:
$${\rho}: B{\Gamma} \to S^1$$
is a bundle. Consider the commutative diagram
$$\begin{CD}
E{\Gamma}{\times}_{\Gamma}{\E}{\G}_{\Gamma} @>>> E{\Gamma}{\times}_{\Gamma}{\E}{\G}_{\Gamma} \\
@V{q_{\Gamma}}VV @VV{p_{\Gamma}}V \\
S^1{\times}{\B}{\G}_{\Gamma} @>>> {\B}{\G}_{\Gamma}
\end{CD}$$
In the applications the spectrum $\bbS$ has an infinite loop space structure. 
Let ${\C}{\F}$ be the homotopy cofiber:
$${\bbS}(BG) \;\xrightarrow{1 - {\phi}_*}\; {\bbS}(BG) \;\to\; {\C}{\F}$$
Since $\text{id}_{BG}$ and $\phi$ are homotopic in $B{\Gamma}$, there is an induced map:
$$f: {\C}{\F} \to {\bbS}(B{\Gamma}).$$
This is the map that `forgets the control' over $S^1$. We will give a homological description of
${\C}{\F}$.

With the above notation:
\begin{itemize}
\item Every simplex $\sigma$ of ${\E}{\G}_{\Gamma}$ defines an equivalence relation on $\bbZ$,
$$m_1 {\sim}_{\sigma}  m_2 \;\Longleftrightarrow\; {\Gamma}_{\sigma}{\cap}(1, m_1)G(1, -m_2) \not=
\emptyset.$$
\item For each $t\in S^1$, $m\in\bbZ$, a map
$$i_m: (EG{\times}\{t\}{\times}\{m\}){\times}_{{\Gamma}_{\sigma}{\cap}G}{\hat{\sigma}} \to
(EG{\times}\{t\}{\times}{\bbZ}){\times}_{{\Gamma}}{\Gamma}{\hat{\sigma}}, \;
[(x, t, m), \hat{\sigma}] \mapsto [(x, t, m), \hat{\sigma}]$$
here $EG{\times}\{t\}{\times}{\bbZ} \subset E{\Gamma}$. 
Notice that the image of $i_m$ depends on the 
choice of $t\in S^1$.
\end{itemize}

\begin{lem}\label{lem-inverse}
Fix $t\in S^1$.
Let $\sigma$ be a simplex of ${\E}{\G}_{\Gamma}$ and $t\in S^1$. Then:
\begin{enumerate}
\item $i_m$ is a monomorphism for each $m\in\bbZ$.
\item $\text{Im}(i_{m_1}) = \text{Im}(i_{m_2})$ if and only if $m_1 {\sim}_{\sigma} m_2$.
\item The images of $i_{m_1}$ and $i_{m_2}$ are either equal or disjoint.
\item The induced map:
$$i = \coprod_{[m_k]}i_{m_k}: \coprod_{[m_k]}\text{Im}(i_k) \to 
(EG{\times}\{t\}{\times}{\bbZ}){\times}_{{\Gamma}}{\Gamma}{\hat{\sigma}}$$
where $m_k$ runs over a complete set of representatives of ${\sim}_{\sigma}$,
is a homeomorphism.
\end{enumerate}
\end{lem}

\begin{proof}
For (1), assume that $i_m([(x, t, m), \hat{\sigma}]) = i_m([(y, t, m), \hat{\sigma}])$. Then,
there is $(g, n)\in\Gamma$, such that
$$((x, t, m), \hat{\sigma}) = ((y, t, m)(g, n), (g, n)^{-1}\hat{\sigma}) = ((y{\alpha}^m(g), 
t, m + n), (g, n)^{-1}\hat{\sigma}).$$
That implies that $(g, n)\in {\Gamma}_{\sigma}$ and $n = 0$. Thus $(g, n) \in {\Gamma}_{\sigma}
{\cap}G$ and $[(x, t, m), \hat{\sigma}] = [(y, t, m), \hat{\sigma}]$. 

Similar calculations show that in $\text{Im}(i_{m_1}) = \text{Im}(i_{m_2})$ then 
$m_1 \sim_{\sigma} m_2$. If we assume that $m_1{\sim}_{\sigma}m_2$, then there is $g\in G$
such that $(g, m_1 - m_2) \in {\Gamma}_{\sigma}$. Then, if $[(x, t, m_1), \hat{\sigma}] \in 
\text{Im}(i_{m_1})$, then
$$[(x, t, m_1), \hat{\sigma}] = [(x, t, m_1), (g, m_1 - m_2)\hat{\sigma}] =
[(x{\alpha}^{m_1}(g), t, m_2), \hat{\sigma}] \in \text{Im}(i_{m_2}),$$
proving (2).

For (3), assume that $m_1$ and $m_2$ are not equivalent and let
$[(x, mt), \hat{\sigma}] \in \text{Im}(i_{m_1}){\cap}\text{Im}(i_{m_2})$.
But by (2), $m {\sim}_{\sigma} m_1$ and $m {\sim}_{\sigma} m_2$. Thus $m_1 {\sim}_{\sigma} m_2$,
contrary to our assumption.

It is immediate that $i$ is surjective. Parts (1)--(3) show that $i$ is a bijection. The inverse
of $i$ is given by:
$$i^{-1} : (EG{\times}\{t\}{\times}{\bbZ}){\times}_{{\Gamma}}{\Gamma}{\hat{\sigma}} \to 
\coprod_{[m_k]}\text{Im}(i_k), \quad 
[(x, t, m), (g, n)\hat{\sigma}] \mapsto [(x{\alpha}^m(g), t, m + n), 
\hat{\sigma}] \in \text{Im}(i_{m+n}),$$
and thus $i$ is a homeomorphism.
\end{proof}

The following result is immediate from the proof.

\begin{lem}\label{lem-induced}
Fix $\sigma\in {\E}{\G}_{\Gamma}$ and $t \in S^1$. Then the map
$${\chi}: EG{\times}_GG(1, m)\hat{\sigma} \to EG{\times}_{\Gamma}{\E}{\G}_{\Gamma},
\; [x, (g, m)\hat{\sigma}] \mapsto [(xg, t, m), \hat{\sigma}]$$
induces a homeomorphism onto $\text{Im}(i_m)$.
\end{lem}

We will study the homotopy colimit structure of the quotient map:
$$q_{\Gamma}: E{\Gamma}{\times}_{\Gamma}{\E}{\G}_{\Gamma} \to {\B}{\G}_{\Gamma}{\times}S^1.$$
Equip $S^1$ with the structure of a simplicial complex 
with three $0$-simplices $v_i$, 
$i = 0, 1, 3$ and three $1$-simplices $e_i = \{v_i, v_{i+1}\}$, where $i$ is taken mod $3$.
The projection map induces a functor:
$${\rho}: \text{cat}(S^1{\times}B{\G}_{\Gamma})^{op} \to \text{cat}(S^1)^{op}.$$
Using Segal's Pushdown Construction we get a homotopy equivalence:
$$\text{hocolim}_{\text{cat}(S^1{\times}B{\G}_{\Gamma})^{op}}
({\bbS}{\circ}\text{bar}(q_{\Gamma})) \simeq 
\text{hocolim}_{\text{cat}(S^1)^{op}}(P_*{\rho}).$$
We will describe explicitely the functor $P_*{\rho}$ on $\text{cat}(S^1)^{op}$:
\par\noindent
For each simplex $t$ of $S^1$, let 
$${\bbS}{\circ}\text{bar}(q_{\Gamma})(t, -): \text{cat}(B{\G}_{\Gamma})^{op} \to 
{\bbS}pectra, \; {\sigma} \mapsto {\bbS}{\circ}\text{bar}(q_{\Gamma})(t, {\sigma})$$
and thus
$$P_*{\rho}(\hat{t}) = \text{hocolim}_{\text{cat}(B{\G}_{\Gamma})^{op}}
({\bbS}{\circ}\text{bar}(q_{\Gamma})(t, -)).$$
So if we fix $t$, the functor $\text{bar}(q_{\Gamma})$ associates to $\sigma$, the space
$$q_{\Gamma}^{-1}(\hat{t}{\times}\hat{\sigma}) = (EG{\times}\{\hat{t}\}{\times}{\bbZ})
{\times}_{\Gamma}{\Gamma}\hat{\sigma} 
\cong (EG{\times}\{\hat{t}\}{\times}{\bbZ})
{\times}_{{\Gamma}_{\sigma}}\hat{\sigma},$$
after subdividing.

The inclusion map $G \to \Gamma$ induces a commutative diagram:
$$\begin{CD}
EG{\times}_G{\E}{\G}_{\Gamma} @>{f}>> 
(EG{\times}\{\hat{t}\}{\times}{\bbZ}){\times}_{\Gamma}{\E}{\G}_{\Gamma} \\
@V{p_G}VV @VV{q_{\Gamma}}V \\
{\E}{\G}_{\Gamma}/G = {\B}{\G}_G @>{u}>> {\hat t}{\times}{\B}{\G}_{\Gamma} = 
{\hat t}{\times}{\E}{\G}_{\Gamma}/{\Gamma}
\end{CD}$$
where $u$ is just the quotient map. Notice that ${\E}{\G}_{\Gamma}$ is a model for
the space of type ${\E}{\G}_G$. Then $u$ induces a functor:
$$U: \text{cat}({\B}{\G}_G)^{op} \to \text{cat}({\B}{\G}_{\Gamma})^{op}$$
Let $[{\sigma}]_{\Gamma}$ denote the $\Gamma$-orbit of the simplex $\sigma$ of
${\E}{\G}_{\Gamma}$. Then the over category:
$$U{\downarrow}[{\sigma}]_{\Gamma} = 
\{[{\gamma}{\tau}]_G: {\tau} \ge {\sigma}, {\gamma}\in{\Gamma}\} = \{[(1, m){\tau}]_G: {\tau} 
\ge {\sigma}, m\in \bbZ\},$$
where $[-]_G$ denotes the $G$-orbit of the simplex.

\begin{lem}\label{lem-claim}
With the above notation,
\begin{enumerate}
\item $\displaystyle{
|U{\downarrow}[{\sigma}]_{\Gamma}| = \coprod_{[m_k]}
|\text{Id}{\downarrow}[(1, m_k){\sigma}]_G|}$.
\item $\displaystyle{
p_G^{-1}(|U{\downarrow}[{\sigma}]_{\Gamma}|)}$ = 
$\displaystyle{\coprod_{[m_k]}
p_G^{-1}(|\text{Id}{\downarrow}[(1, m_k){\sigma}]_G|)} =$ \\
$\displaystyle{\coprod_{[m_k]}
EG{\times}_GG(|\text{Id}{\downarrow}(1, m_k){\sigma}|) \cong 
\coprod_{[m_k]}
EG{\times}_{{\Gamma}_{\sigma}{\cap}G}|\text{Id}{\downarrow}(1, m_k){\sigma}|}$.
\end{enumerate}
where $m_k$ runs over a complete set of representatives of $\sim_{\sigma}$.
\end{lem}

\begin{proof} 
Notice that:
$$[ {\gamma}_1{\sigma}]_G = [{\gamma}_2{\sigma}]_G \;\Longleftrightarrow\;
{\Gamma}_{\sigma}{\cap}{\gamma}_1G{\gamma}_2^{-1} \not= \emptyset$$
In this case, each orbit admits a representative of the form $(1, m){\tau}$, with 
$m\in\bbZ$,
$$[(1, m_1){\sigma}]_G = [(1, m_2){\sigma}]_G \;\Longleftrightarrow\;
{\Gamma}_{\sigma}{\cap}(1, m_1)G(1, -m_2) \not= \emptyset \;\Longleftrightarrow\; 
m_1{\sim}_{\sigma}m_2,$$
which implies that $|U{\downarrow}[{\sigma}]_{\Gamma}|$ is equal to the union on the right
side. Since we used the double subdivision of the original object, it follows that $\sigma$
and ${\gamma}{\sigma}$ are equal or there is no chain of simplices containing both of them. 
That proves Part (1).

The first equality for Part (2) follows from (1), and the last homeomorphism is standard.
\end{proof}

Now we have all the basic tools to apply Segal's Theorem on the functor $U$.

\begin{lem}\label{lem-basic}
Assume that $G$ satisfies the $\bbS$-IC. Then
for each $t \in \text{cat}(S^1)^{op}$, there is a natural homotopy equivalence:
$$P_*{\rho}(t) = \text{hocolim}_{\text{cat}(B{\G}_{\Gamma})^{op}}
({\bbS}{\circ}\text{bar}(q_{\Gamma})(t, -)) \simeq {\bbS}(BG).$$
\end{lem}

\begin{proof}
We will show that
$$\text{hocolim}_{\text{cat}({\B}{\G})_G^{op}}({\bbS}{\circ}\text{bar}(p_G)) \simeq
\text{hocolim}_{\text{cat}(B{\G}_{\Gamma})^{op}}
({\bbS}{\circ}\text{bar}(q_{\Gamma})(t, -)).$$
Since $G$ satisfies the $\bbS$-IC, the result will follow.

Segal's Theorem implies that
$$\text{hocolim}_{\text{cat}({\B}{\G})_G^{op}}({\bbS}{\circ}\text{bar}(p_G)) \simeq
\text{hocolim}_{\text{cat}({\B}{\G})_{\Gamma}^{op}}(h)$$
where, for each simplex $\sigma$ of ${\E}{\G}_{\Gamma}$,
$$h([{\sigma}]_{\Gamma}) = 
\text{hocolim}_{\text{cat}(|U{\downarrow}[{\sigma}]_{\Gamma}|)^{op}}({\bbS}{\circ}
\text{bar}(p_G|))$$  
where $[{\sigma}]_{\Gamma}$ is the $\Gamma$-orbit of $\sigma$ and 
$$p_G|: (p_G)^{-1}(|U{\downarrow}[{\sigma}]_{\Gamma}|) \to |U{\downarrow}[{\sigma}]_{\Gamma}|$$
is the restriction of $p_G$. Using Lemma \ref{lem-claim},
$$\begin{array}{llll}
h([{\sigma}]_{\Gamma}) & = &
\text{hocolim}_{\text{cat}(|U{\downarrow}[{\sigma}]_{\Gamma}|)^{op}}({\bbS}{\circ}
\text{bar}(p_G|)) & \\[.5ex]
&\simeq & \displaystyle{\bigvee_{[m_k]}
\text{hocolim}_{\text{cat}(|\text{Id}{\downarrow}[(1, m_k){\sigma}]_G|)^{op}}
({\bbS}{\circ}
\text{bar}(p_G|))}, & \text{from Lemma \ref{lem-claim}, Part (1)}\\[.5ex]
& \simeq & \displaystyle{\bigvee_{[m_k]}
\text{hocolim}_{\text{cat}(|[(1, m_k){\sigma}]_G|)^{op}}
({\bbS}{\circ}
\text{bar}(p_G|))}, & \text{because $|\text{Id}{\downarrow}[(1, m_k){\sigma}]_{\Gamma}| \simeq
[(1, m_k){\sigma}]_G$}\\[.5ex]
& \simeq & \displaystyle{\bigvee_{[m_k]}{\bbS}(p_G^{-1}([(1, m_k)\hat{\sigma}]_G))}, & 
\text{from the definition} \\[.5ex]
& \simeq & \displaystyle{{\bbS}\left(\coprod_{[m_k]}p_G^{-1}([(1, m_k)\hat{\sigma}]_G)\right)} & 
\\[.5ex]
& = & \displaystyle{{\bbS}\left(\coprod_{[m_k]}EG{\times}_GG(1, m_k)\hat{\sigma}\right)} & \\[.5ex]
& \simeq & \displaystyle{{\bbS}\left(\coprod_{[m_k]}\text{Im}(i_k) \right)}, & \text{from Lemma 
\ref{lem-induced}} \\[.5ex]
& \simeq & {\bbS}((EG{\times}\{\hat{t}\}{\times}{\bbZ})
{\times}_{\Gamma}{\Gamma}\hat{\sigma}), 
& \text{from Lemma \ref{lem-inverse}} \\[.5ex]
& = & {\bbS}{\circ}\text{bar}(q_{\Gamma})(t, 
[{\sigma}]_{\Gamma}) &
\end{array}
$$ 
(the wedge is the coproduct in the category of spectra).
Therefore 
$${\bbS}(BG) \simeq 
\text{hocolim}_{\text{cat}({\B}{\G})_G^{op}}({\bbS}{\circ}\text{bar}(p_G)) \simeq
\text{hocolim}_{\text{cat}({\B}{\G})_{\Gamma}^{op}}(h) \simeq 
\text{hocolim}_{\text{cat}({\B}{\G})_{\Gamma}^{op}}(
{\bbS}{\circ}\text{bar}(q_{\Gamma})(t, -))$$
completing the proof.
\end{proof}

\begin{prop}\label{prop-circle}
If the ${\bbS}$-IC holds for $G$, then the homotopy cofiber of $1 - {\phi}_*$ is 
${\bbH}.(S^1{\times}B{\G}_{\Gamma}, {\bbS}(q_{\Gamma}))$. Thus there is a long exact sequence:
$$\dots \to {\pi}_i({\bbS}(BG)) \xrightarrow{1 - {\phi}_*} {\pi}_i({\bbS}(BG)) \to
H_i(S^1{\times}B{\G}_{\Gamma}, {\bbS}(q_{\Gamma})) \to {\pi}_{i-1}({\bbS}(BG)) \to \dots$$
\end{prop}

\begin{proof}
We use the fact that $q_{\Gamma}$ has a homotopy colimit structure and
$$H_i(S^1{\times}B{\G}_{\Gamma}, {\bbS}(q_{\Gamma})) \cong
{\pi}_i(\text{hocolim}_{\text{cat}(S^1{\times}B{\G}_{\Gamma})^{op}}
({\bbS}{\circ}\text{bar}(q_{\Gamma})))$$
Using the multiplicative properties of homotopy colimits and Lemma \ref{lem-basic}, 
$$H_i(S^1{\times}B{\G}_{\Gamma}, {\bbS}(q_{\Gamma})) = 
{\pi}_i(\text{hocolim}_{\text{cat}(S^1)^{op}}({\bbS}(BG)))$$
The result follows as in Section 3 in \cite{mupr:wa}.
\end{proof}

The commutative diagram
$$\begin{CD}
E{\Gamma}{\times}_{\Gamma}{\E}{\G}_{\Gamma} @>>> E{\Gamma}{\times}_{\Gamma}{\E}{\G}_{\Gamma} 
@>>> B{\Gamma} \\
@V{q_{\Gamma}}VV @VV{p_{\Gamma}}V @VVV \\
S^1{\times}{\B}{\G}_{\Gamma} @>>> {\B}{\G}_{\Gamma} @>>> *
\end{CD}$$
induces a map of spectra
$$f: {\bbH}.(S^1{\times}B{\G}_{\Gamma}, {\bbS}(q_{\Gamma})) \to {\bbS}(B{\Gamma}).$$
 
\begin{cor}\label{cor-com}
The following diagram commutes, up to homotopy:
$$\begin{CD}
{\bbS}(BG) @>>> {\bbH}.(S^1{\times}B{\G}_{\Gamma}, {\bbS}(q_{\Gamma})) \\
@V{\text{id}}VV @VV{f}V \\
{\bbS}(BG) @>>> {\bbS}(B{\Gamma})
\end{CD}$$
\end{cor}

\begin{proof}
This is because the inclusion induced map ${\bbS}(BG) \to {\bbS}(B{\Gamma})$ 
factors through the spectrum
${\bbH}.(S^1{\times}B{\G}_{\Gamma}, {\bbS}(q_{\Gamma}))$ because of the commutative diagram
$$\begin{CD}
BG @<<< EG{\times}_G{\E}{\G}_{\Gamma} @>>>
E{\Gamma}{\times}_{\Gamma}{\E}{\G}_{\Gamma} @>>> E{\Gamma}{\times}_{\Gamma}{\E}{\G}_{\Gamma} 
@>>> B{\Gamma} \\
@VVV @V{p_G}VV @V{q_{\Gamma}}VV @VV{p_{\Gamma}}V @VVV \\
* @<<< {\B}{\G}_G @>>> S^1{\times}{\B}{\G}_{\Gamma} @>>> {\B}{\G}_{\Gamma} @>>> *
\end{CD}$$
Notice that the map in the second square is not natural because it depends on the choice of 
an element of $S^1$ but
the composition induced by the third square is not affected by that choice. 
\end{proof}

We set up the notation for the Bundle version of Proposition \ref{prop-circle}. Start
with a commutative diagram:
$$\begin{CD}
\bar{\bar{Y}} @>>> \bar{\bar{Y}} @>>> \bar{Y} @>>> Y \\
@V{\bar{\bar{\rho}}}VV @V{\bar{\bar{\rho}}}VV @V{\bar{\rho}}VV @VV{\rho}V \\
EG{\times}_G{\E}{\G}_{\Gamma} @>>>
EG{\times}_G{\E}{\G}_{\Gamma} @>>> E{\Gamma}{\times}_{\Gamma}{\E}{\G}_{\Gamma} @>>> B{\Gamma} \\
@V{q_G}VV @V{p_G}VV @V{p_{\Gamma}}VV @VVV \\
S^1{\times}{\B}{\G}_G @>>> {\B}{\G}_G @>>> {\B}{\G}_{\Gamma} @>>> *
\end{CD}$$
where $\rho$ is a bundle and the top diagrams are pull-back diagrams. Then $\bar{\bar{\rho}}$ 
is also a bundle. The proofs of Lemma \ref{lem-basic} and Proposition \ref{prop-circle} 
also work in this case. 

\begin{thm}\label{thm-circle}
Assume that the Bundle $\bbS$-IC holds for $G$. Then there is an exact sequence:
$$\dots \to {\pi}_i({\bbS}(\bar{\bar{Y}})) \xrightarrow{1 - {\phi}_*} 
{\pi}_i({\bbS}(\bar{\bar{Y}})) \to
H_i(S^1{\times}B{\G}_{\Gamma}, {\bbS}(\bar{\rho}{\circ}q_{\Gamma})) 
\to {\pi}_{i-1}({\bbS}(\bar{\bar{Y}})) \to \dots$$
where $\phi: \bar{\bar{Y}} \to \bar{\bar{Y}}$ is the homeomorphism induced by $\alpha$.
\end{thm}

We specialize to the case ${\bbS} = {\bbK}_R$, the $K$-theory spectrum with coefficients
in a ring $R$. Then Corollary \ref{cor-com} implies that there is a commutative diagram 
of exact sequences for each $i$:
$$\xymatrix{
& & H_i(S^1{\times}{\B}{\G}_{\Gamma}, {\bbK}_R(q_{\Gamma})) \ar[dr] \ar[dd]^f& & \\
K_i(RG) \ar[r]^{1 - {\phi}_*} & K_i(RG) \ar[dr] \ar[ur] & &  
K_{i-1}(RG) \ar[r]^{1 - {\phi}_*}
& K_{i-1}(RG)\\
& & K_i(R{\Gamma})  \ar[ur] & &
}$$

\begin{prop}\label{prop-iso}
Let $RG$ be a regular coherent ring. Then
$$f: H_i(S^1{\times}{\B}{\G}_{\Gamma}, {\bbK}_R(q_{\Gamma})) \to K_i(R{\Gamma})$$
is an isomorphism for all $i\in \bbZ$.
\end{prop}

\begin{proof}
This follows from the analogue of the Bass--Heller--Swan formula for semidirect products
(\cite{wa:wh}, \cite{wa:al}). The top sequence is exact by Proposition \ref{prop-circle}. The
bottom sequence is exact because the assumption on $R{\Gamma}$ implies that the Nil-groups
vanish.
\end{proof}

>From the definition of $f$ we have that $f$ factors as:
$$f: H_i(S^1{\times}{\B}{\G}_{\Gamma}, {\bbK}_R(q_{\Gamma})) \xrightarrow{p}
H_i({\B}{\G}_{\Gamma}, {\bbK}_R(p_{\Gamma})) \xrightarrow{A} K_i(R{\Gamma})$$
where $p$ is induced by the projection to the second coordinate and $A$ is the
assembly map.

\begin{prop}\label{prop-epi}
Let $RG$ be a regular coherent ring. Then
$$A: H_i({\B}{\G}_{\Gamma}, {\bbK}_R(p_{\Gamma})) \to K_i(R{\Gamma})$$
is an epimorphism.
\end{prop}

\begin{proof}
It follows from Proposition \ref{prop-iso}.
\end{proof}

\section{A Special Case}\label{sec-trivial}

Let ${\Gamma} = G{\times}{\bbZ}$. Then a model for $E{\Gamma}$ can be chosen to be 
$EG{\times}{\bbR}$, where $EG$ is any model for the classifying space of $G$. In this case,
Proposition \ref{prop-circle} has a much simpler interpretation. The assumption for this section
is that:

\noindent
{\bf Assumption (IC)}: The ${\bbS}$-IC holds for groups of the form $G_0{\times}{\bbZ}$ where
$G_0$ is a virtually infinite cyclic subgroup of $G$.

\begin{rem}
For ${\bbS}$ the pseudoisotopy spectrum the result follows from \cite{fajo:is} because the
$G_0{\times}{\bbZ}$ is virtually abelian. For the $\bbK$-theory spectrum it is not known
if the $\bbK$-IC holds for such groups. Partial results in this direction are in \cite{barfajore}
and \cite{fali}.
\end{rem} 

\begin{prop}\label{prop-circle-trivial}
Let $\Gamma = G{\times}{\bbZ}$. Assume that the $\bbS$-IC holds for $G$. 
Then, in the notation of the last section, there is a homotopy equivalence:
$${\bbH}.(S^1{\times}B{\G}_{\Gamma}, {\bbS}(q_{\Gamma})) \simeq 
{\bbS}(BG){\times}{\Omega}^{-1}({\bbS}(BG)).$$
\end{prop}

The projection map ${\Gamma} \to G$ induces a commutative diagram
$$\begin{CD}
(EG{\times}{\bbR}){\times}_{\Gamma}{\E}{\G}_{\Gamma} @>>> (EG{\times}{\bbR}){\times}_{\Gamma}
{\E}{\G}_G \\
@V{p_{\Gamma}}VV @VV{q}V \\
{\B}{\G}_{\Gamma} @>{u}>> {\B}{\G}_G
\end{CD}$$
where $\Gamma$ acts on ${\E}{\G}_G$ through the projection ${\Gamma} \to G$.

\begin{lem}\label{lem-inv}
Let $[x]$ be a point in ${\B}{\G}_G$. Then $q^{-1}([x]) \cong p_G^{-1}([x]){\times}S^1$.
\end{lem}

\begin{proof}
Let $[x]{\in}{\B}{\G}_G$, with $x$ the $G$-orbit of a point $x{\in}{\E}{\G}_G$. 
Then a point $[(e, r), y]{\in} (EG{\times}{\bbR}){\times}_{\Gamma}
{\E}{\G}_G$ is in $q^{-1}(x)$, if there is $g{\in}G$ such that $gy = x$. Define a map
$${\alpha}: q^{-1}([x]) \to p_G^{-1}([x]){\times}S^1, \;
[(e, r), y] \mapsto ([e, y], [r])$$
where $[r]$ is the $\bbZ$-orbit or $r$. The map is a well-defined homeomorphism.
\end{proof}

The map $u$ induces
a functor $U: \text{cat}({\B}{\G}_{\Gamma})^{op} \to \text{cat}({\B}{\G}_G)^{op}$.
Let $\sigma$ be a simplex ${\B}{\G}_G$ with isotropy group $G_0$. 
We will describe the category
$U{\downarrow}{\sigma}$. Let $[\sigma]_G$ be the $G$-orbit of a simplex $\sigma$ of
${\E}{\G}_G$ with isotropy group $G_0 < G$.
A simplex $[\bar{\sigma}]_{\Gamma}$ of 
${\B}{\G}_{\Gamma}$ will be an object of $U{\downarrow}{\sigma}$ if
$u([\bar{\sigma}]_{\Gamma}) \ge [{\sigma}]_G$.
The homotopy invariant functor $\bbS$ induces a functor
$${\bbS}[S^1]: {\bf Top} \to {\bbS}pectra, \; {\bbS}[S^1](X) = {\bbS}(X{\times}S^1)$$
which is also homotopy invariant.

\begin{prop}\label{prop-trivial}
Assume that $G$ satisfies Assumption (IC).
With the above notation, there is a homotopy equivalence of spectra:
$$\text{hocolim}_{\text{cat}({\B}{\G}_{\Gamma})^{op}}({\bbS}{\circ}\text{bar}(p_{\Gamma})) \simeq
\text{hocolim}_{\text{cat}({\B}{\G}_G)^{op}}({\bbS}[S^1]{\circ}\text{bar}(p_G)).$$
\end{prop}

\begin{proof}
We will use Segal's Pushdown Theorem for comparing the two homotopy colimits. The map $u$ induces
a functor $U: \text{cat}({\B}{\G}_{\Gamma})^{op} \to \text{cat}({\B}{\G}_G)^{op}$. Let
$[{\sigma}]_G$ be a simplex in ${\B}{\G}_G$, with isotropy group $G_0$ (defined up to conjugation).
We choose the classifying spaces of $G_0$ and $G_0{\times}{\bbZ}$:
\begin{enumerate}
\item Since $G_0$ is in ${\G}_G$, 
${\B}{\G}_{G_0}$ can be chosen to be a representative $\sigma$ and
${\E}{\G}_{G_0} = |{\sigma}|$. 
\item Since $G_0{\times}{\bbZ} < {\Gamma}$, we choose ${\E}{\G}_{G_0{\times}{\bbZ}}$
to be the subcomplex of ${\E}{\G}_{\Gamma}$ consisting of all simplices with isotropy 
group a subgroup of a $G_0{\times}{\bbZ}$-conjugate of $G_0{\times}{\bbZ}$ i.e. 
of the form ${\gamma}H{\gamma}^{-1}$, with ${\gamma}\in G_0{\times}{\bbZ}$ and 
$H < G_0{\times}{\bbZ}$.
Set ${\B}{\G}_{G_0{\times}{\bbZ}} = {\E}{\G}_{G_0{\times}{\bbZ}}/G_0{\times}{\bbZ}$.
\end{enumerate}
With the above choices, we get a commutative diagram:
$$\begin{CD}
(EG{\times}{\bbR}){\times}_{G_0{\times}{\bbZ}}{\E}{\G}_{G_0{\times}{\bbZ}} @>>> 
(EG{\times}{\bbR}){\times}_{G_0{\times}{\bbZ}}|{\sigma}| @>{\simeq}>> B(G_0{\times}{\bbZ}) \\
@V{p_{G_0{\times}{\bbZ}}}VV @V{q_{G_0{\times}{\bbZ}}}VV @VVV \\
{\B}{\G}_{G_0{\times}{\bbZ}} @>{u_0}>> |[{\sigma}]_G| @>>> *
\end{CD}$$
Since the ${\bbS}$-IC holds for $G_0{\times}{\bbZ}$, the induced map
$$\text{hocolim}_{\text{cat}({\B}{\G}_{G_0{\times}{\bbZ}})^{op}}({\bbS}{\circ}
\text{bar}(p_{{G_0}{\times}{\bbZ}})) \to
\text{hocolim}_*({\bbS}{\circ}*) = {\bbS}(B(G_0{\times}{\bbZ}))$$
is a homotopy equivalence. Since in the right square the horizontal maps are homotopy equivalences,
the induced maps
$$\text{hocolim}_{\text{cat}({\B}{\G}_{G_0{\times}{\bbZ}})^{op}}({\bbS}{\circ}
\text{bar}(p_{G_0{\times}{\bbZ}})) 
\to \text{hocolim}_{\text{cat}([{\sigma}]_G)^{op}}({\bbS}{\circ}
\text{bar}(q_{G_0{\times}{\bbZ}})) 
\simeq {\bbS}(B(G_0{\times}{\bbZ})) \quad (*)$$
are also homotopy equivalences. The map $u_0$ induces a functor $U_0$ between the corresponding 
categories. We will apply Segal's Pushdown Theorem on the functor $U_0$. For each 
$[{\tau}]_G < [{\sigma}]_G$, 
define the following map, which is the restriction of $p_{G_0{\times}{\bbZ}}$:
$$p_{\tau}: (EG{\times}{\bbR}){\times}_{G_0{\times}\bbZ}p_{G_0{\times}{\bbZ}}^{-1}
(|U_0{\downarrow}[{\tau}]_G|)
\to |U_0{\downarrow}[{\tau}]_G|$$
Define a functor 
$$h_{\sigma}:  \text{cat}({[\sigma}]_G)^{op} \to {\bbS}pectra, \;
h({\tau}) = \text{hocolim}_{U_0{\downarrow}[{\tau}]_G}({\bbS}{\circ}\text{bar}(p_{\tau})).$$
Then Segal's Pushdown Construction implies that
$$\text{hocolim}_{\text{cat}({\B}{\G}_{G_0{\times}{\bbZ}})^{op}}({\bbS}{\circ}
\text{bar}(p_{G_0{\times}{\bbZ}}))
\simeq  \text{hocolim}_{\text{cat}([{\sigma}]_G)^{op}}h.$$
But the category $\text{cat}([{\sigma}]_G)^{op}$ 
has a unique minimal element, namely $[\sigma]_G$. Thus
$$\text{hocolim}_{\text{cat}({\B}{\G}_{G_0{\times}{\bbZ}})^{op}}({\bbS}{\circ}
\text{bar}(p_{G_0{\times}{\bbZ}}))
\simeq  \text{hocolim}_{\text{cat}([{\sigma}]_G)^{op}}h \simeq 
\text{hocolim}_{U_0{\downarrow}[{\sigma}]_G}({\bbS}{\circ}\text{bar}(q_{G_0{\times}{\bbZ}})).$$
Each simplex of $U{\downarrow}[{\sigma}]_G$ has isotropy group that it is contained  in 
$G_0{\times}{\bbZ}$.
Thus, the natural map induces a commutative diagram
$$\begin{CD}
(EG{\times}{\bbR}){\times}_{G_0{\times}
{\bbZ}}p_{G_0{\times}{\bbZ}}^{-1}(|U_0{\downarrow}[{\sigma}]_G|) @>>>
(EG{\times}{\bbR}){\times}_{\Gamma}p_{\Gamma}^{-1}(|U{\downarrow}[{\sigma}]_G|) \\
@V{p_{G_0{\times}{\bbZ}}|}VV @VV{p_{\Gamma}|}V \\
|U_0{\downarrow}[{\sigma}]_G| @>>> |U{\downarrow}[{\sigma}]_G|
\end{CD}$$
and the horizontal maps are homeomorphisms. Therefore,
$$\text{hocolim}_{U_0{\downarrow}[{\sigma}]_G}({\bbS}{\circ}\text{bar}(p_{G_0{\times}{\bbZ}})) 
\simeq 
\text{hocolim}_{U{\downarrow}[{\sigma}]_G}({\bbS}{\circ}\text{bar}(p_{\Gamma})).$$
Combining with (*),
$${\bbS}(B(G_0{\times}{\bbZ})) \simeq 
\text{hocolim}_{U{\downarrow}[{\sigma}]_G}({\bbS}{\circ}\text{bar}(p_{G_0{\times}{\bbZ}}))$$
which implies that
$${\bbS}(B(G_0{\times}{\bbZ})) \simeq 
\text{hocolim}_{U{\downarrow}[{\sigma}]_G}({\bbS}{\circ}\text{bar}(p_{G_0{\times}{\bbZ}})) 
\quad (**)$$
and the homotopy equivalence is natural for $[\sigma]_G$ a simplex in ${\B}{\G}_G$.

Now we apply Segal's Pushdown Construction to the functor $U$. Then
$$\text{hocolim}_{\text{cat}({\B}{\G}_{\Gamma})^{op}}({\s}{\circ}\text{bar}(p_{\Gamma})) \simeq
\text{hocolim}_{\text{cat}({\B}{\G}_G)^{op}}({\chi})$$
where, for each simplex $[\sigma]_G$ of ${\B}{\G}_G$, 
$${\chi}({\sigma}) = \text{hocolim}_{U{\downarrow}[{\sigma}]_G}({\bbS}{\circ}
\text{bar}(p_{G_0{\times}{\bbZ}})).$$
>From (**), we have that 
 $$\text{hocolim}_{\text{cat}({\B}{\G}_{\Gamma})^{op}}({\bbS}{\circ}\text{bar}(p_{\Gamma})) \simeq
\text{hocolim}_{\text{cat}({\B}{\G}_G)^{op}}({\psi})$$
where ${\psi}([{\sigma}]_G) = 
{\bbS}(B(G_{\sigma}{\times}{\bbZ}))$. The result follows from the
definition and the homotopy invariance of the functor $\bbS$.
\end{proof}

We now specialize to the case that $\bbS = {\bbK}_R$ 
is the $K$-theory spectrum with coefficients in
a ring $R$. Then Proposition \ref{prop-trivial} provides an alternative description of
the homology term in the IC, corresponding to the Bass-Heller-Swan splitting of $K$-theory
(\cite{ba:al}, \cite{gr}, \cite{ra:lo}).
For the definition of the Nil-spectrum functor we use ideas from the Bass-Heller-Swan Formula:
$${\bbN}il_R(-) = \text{Cofiber}[{\bbK}_R(-) \to {\bbK}_{R[t]}(-)]$$
where $R[t]$ is the polynomial ring. The Bass-Heller-Swan Formula implies that 
there are natural homotopy equivalences:
$$\begin{array}{rll}
{\bbK}_{R[t]}(-) & \simeq &  {\bbK}_R(-) {\times} {\bbN}il_R(-) \\
{\bbK}_{R[t^{-1}]}(-) & \simeq &  {\bbK}_R(-) {\times} {\bbN}il_R(-) \\
{\bbK}_{R[t,t^{-1}]}(-) & \simeq &  {\bbK}_R(-) {\times} {\Omega}^{-1} {\bbK}_R(-) 
{\times} {\bbN}il_R(-) {\times} {\bbN}il_R(-) 
\end{array}$$

The following result is the homological analogue of the Bass-Heller-Swan Formula.

\begin{lem}[Bass--Heller--Swan Formula]\label{lem-bhs}
Let $G$ be a group that satisfies Assumption (IC). Then
there is a homotopy equivalence:
$$\begin{array}{ll}
&{\bbH}.({\B}{\G}_{\Gamma}, {\bbK}_R(p_{\Gamma})) \simeq \\
&{\bbH}.({\B}{\G}_G, {\bbK}_R(p_G)) {\times}
{\bbH}.({\B}{\G}_G, {\Omega}^{-1}{\bbK}_R(p_G)) 
{\times}
{\bbH}.({\B}{\G}_G, {\bbN}il_R(p_G)) {\times}
{\bbH}.({\B}{\G}_R, {\bbN}il_R(p_G)).
\end{array}$$
\end{lem}

\begin{proof}
Proposition \ref{prop-trivial} shows that 
$${\bbH}.({\B}{\G}_{\Gamma}, {\bbK}_R(p_{\Gamma})) \simeq 
{\bbH}.({\B}{\G}_G, {\bbK}_{R[t,t^{-1}]}(p_G)).$$
The Bass-Heller-Swan Formula applies to the coefficient spectrum. The result follows.
\end{proof}

The homotopy groups of the homology spectra can be computed using spectral sequences. Using
this description, we immediately have that:
$${\bbH}.({\B}{\G}_G, {\Omega}^{-1}{\bbK}_R(p_G)) \simeq {\Omega}^{-1}{\bbH}.({\B}{\G}_G,
{\bbK}_R(p_G)).$$

\begin{cor}\label{cor-des}
Let $G$ be a group as in Lemma \ref{lem-bhs}. Then
there is a homotopy equivalence:
$$\begin{array}{ll}
&{\bbH}.({\B}{\G}_{\Gamma}, {\bbK}_R(p_{\Gamma})) \simeq \\
&{\bbH}.({\B}{\G}_G, {\bbK}_R(p_G)) {\times}
{\Omega}^{-1}{\bbH}.({\B}{\G}_G, {\bbK}_R(p_G)) 
{\times}
{\bbH}.({\B}{\G}_G, {\bbN}il_R(p_G)) {\times}
{\bbH}.({\B}{\G}_R, {\bbN}il_R(p_G)).
\end{array}$$
\end{cor}

The following summarizes certain immediate consequences of the calculations above.

\begin{cor}\label{cor-cons}
Let $G$ be a torsion free group that satisfies the ${\bbK}_R$-IC and
Assumption (IC). Let ${\Gamma} = G{\times}{\bbZ}$. 
If $R$ is a regular coherent ring, then there are homotopy equivalences:
$${\bbK}_R(BG){\times}{\Omega}^{-1}{\bbK}_R(BG) \simeq
{\bbH}.(S^1{\times}{\B}{\G}_{\Gamma}, {\bbK}_R{\circ}\text{bar}(q'_{\Gamma})) \simeq
{\bbH}.({\B}{\G}_{\Gamma}, {\bbK}_R{\circ}\text{bar}(p_{\Gamma})).$$
Furthermore, the commutative diagram
$$\begin{CD}
E{\Gamma}{\times}{\E}{\G}_{\Gamma} @>>> E{\Gamma}{\times}{\E}{\G}_{\Gamma} \\
@V{q'_{\Gamma}}VV @VV{p_{\Gamma}}V \\
S^1{\times}{\B}{\G}_{\Gamma} @>>>  {\B}{\G}_{\Gamma}
\end{CD}$$
induces the second homotopy equivalence.
\end{cor}

\begin{proof}
Since $G$ is torsion free, the virtually cyclic subgroups of $G$ are infinite cyclic. But, since
$R$ is regular coherent the Nil-groups of $R$ and $R[{\bbZ}]$ vanish (\cite{quillen},
\cite{wa:al}). Therefore the
Nil-spectrum of the corresponding spaces is contractible. Lemma \ref{lem-bhs} 
and Corollary \ref{cor-des} imply that
$$\begin{array}{ll}
&{\bbH}.({\B}{\G}_{\Gamma}, {\bbK}_R(p_{\Gamma})) \simeq 
{\bbH}.({\B}{\G}_G, {\bbK}_R(p_G)) {\times}
{\bbH}.({\B}{\G}_G, {\Omega}^{-1}{\bbK}_R(p_G)) \simeq \\
&{\bbH}.({\B}{\G}_G, {\bbK}_R(p_G)) {\times}
{\Omega}^{-1}{\bbH}.({\B}{\G}_G, {\bbK}_R(p_G)) \simeq {\bbK}_R(BG) {\times}
{\Omega}^{-1}{\bbK}_R(BG).
\end{array}$$
The result follows from Proposition \ref{prop-circle-trivial}.

For the second part of the Corollary, let $\phi$ be the map induced from
the commutative diagram. Notice that there are commutative diagrams:
$$\begin{CD}
EG{\times}_G{\E}{\G}_{\Gamma} @>>> (EG{\times}{\bbZ}t){\times}_{\Gamma}{\E}{\G}_{\Gamma} @>>>
E{\Gamma}{\times}{\E}{\G}_{\Gamma} @>>> 
E{\Gamma}{\times}{\E}{\G}_{\Gamma} \\
@V{p_G}VV @V{q'_{\Gamma}|}VV @V{q'_{\Gamma}}VV @V{p_{\Gamma}}VV \\
{\B}{\G}_G @>>>
\{t\}{\times}{\B}{\G}_{\Gamma} @>>> S^1{\times}{\B}{\G}_{\Gamma} @>>>  {\B}{\G}_{\Gamma}
\end{CD}$$
The composition map induced by the commutative diagrams induce the inclusion
$${\bbK}_R(BG) \to {\bbH}.({\B}{\G}_{\Gamma}, {\bbK}_R(p_{\Gamma})).$$
The map induced from the first two diagrams induce the map
$${\bbK}_R(BG) \to {\bbH}.(S^1{\times}{\B}{\G}_{\Gamma}, {\bbK}_R(q'_{\Gamma})).$$
Thus the diagram commutes:
$$\begin{CD}
{\bbK}_R(BG) @>>> {\bbH}.({\B}{\G}_{\Gamma}, {\bbK}_R(p_{\Gamma})) \\
@| @VV{\phi}V \\
{\bbK}_R(BG) @>>> {\bbH}.(S^1{\times}{\B}{\G}_{\Gamma}, {\bbK}_R(q'_{\Gamma}))
\end{CD}.$$
The result follows.
\end{proof}

\begin{thm}\label{thm-trivial-regular}
Let $RG$ be a regular coherent 
ring, the ${\bbK}_R$-IC holds for $G$, and $G$ satisfies Assumption (IC). 
Then the ${\bbK}_R$-IC holds
for $\Gamma$. 
\end{thm}

\begin{proof}
That follows from the Bass-Heller-Swan Formula for ${\bbK}_R$(\cite{quillen}):
$${\bbK}_R(B{\Gamma}) \simeq {\bbK}_R(BG){\times}{\Omega}^{-1}{\bbK}_R(BG)$$
because $RG$ is regular coherent. The result follows from Corollary \ref{cor-cons} because
$G$ must be torsion free for $RG$ to have finite cohomological dimension. 
\end{proof}

For the Bundle ${\bbK}_R$-IC the regularity of $R{\Gamma}$ is not needed.

\begin{prop}\label{prop-trivial-bundle}
Let $G$ a group that satisfies the bundle ${\bbK}_R$-IC and the bundle version
of Assumption (IC). 
Then the bundle ${\bbK}_R$-IC holds
for $\Gamma$. 
\end{prop}

\begin{proof}
The proof follows from Corollary \ref{cor-bun}.
\end{proof}

Another application of the methods coming from the Bass-Heller-Swan Formula is the
IC for the  ${\bbN}il_R$-spectrum.

\begin{thm}\label{thm-nil}
Let $G$ be a group that satisfies the Bundle ${\bbK}_R$-IC and the bundle version of
Assumption (IC). 
Then $G$ satisfies the 
Bundle ${\bbN}il_R$-IC.
\end{thm}

\begin{proof}
Again we will give the proof when the bundle map is the identity. The general case
follows similarly.
For $F < G$ a virtually cyclic subgroup, its inverse image under the projection map
${\Gamma} \to G$ are of the form $F{\times}{\bbZ}$.
Since $G$ and the inverse image of its virtually cyclic subgroups satisfy the Bundle 
${\bbK}_R$-IC,  $\Gamma$ does too
(Theorem \ref{thm-projection}). Since ${\bbK}_R$-IC holds for $G$, 
Proposition \ref{prop-circle} and Lemma \ref{lem-bhs} imply that
$$\begin{CD}
{\bbH}.(S^1{\times}{\B}{\G}_{\Gamma}, {\bbK}_R(q_{\Gamma})){\times}
{\bbH}.(B{\Gamma}, {\bbN}il_R(p_{\Gamma}))
{\times} {\bbH}.(B{\Gamma}, {\bbN}il_R(p_{\Gamma})) 
@>{\simeq}>> {\bbH}.(B{\Gamma}, {\bbK}_R(p_{\Gamma})) \\
@VVV @VVV \\
{\bbH}.(S^1{\times}{\B}{\G}_{\Gamma}, {\bbK}_R(q_{\Gamma}))
{\times}{\bbN}il_R(B{\Gamma}) {\times} {\bbN}il_R(B{\Gamma}) 
 @>{\simeq}>>  {\bbK}_R(B{\Gamma})
\end{CD}
$$
where the bottom horizontal map is a homotopy equivalence by the Bass-Heller-Swan Splitting.
The vertical maps are induced by the assembly map and the right one is a homotopy equivalence
because the ${\bbK}_R$-IC holds for $\Gamma$. That implies that the left map, which is the
product of the identity and two assembly maps, is an isomorphism. Thus, the assembly map
$${\bbH}.(B{\Gamma}, {\bbN}il_R(p_{\Gamma})) \to {\bbN}il_R(B{\Gamma}) $$
is a homotopy equivalence.
\end{proof}

\begin{rem}
In Theorem \ref{thm-nil} we can not remove the bundle assumption even if only the ${\bbN}il_R$-IC
is to be proved.
\end{rem}

\section{Controlled Groups over the Interval}\label{sec-interval}

Let ${\Gamma} = G_1*_{G_0}G_2$, where $H$ is a subgroup of $G_1{\cap}G_2$. Let $BG_i$, 
$i = 0, 1, 2$, be classifying spaces for the corresponding groups with $BG_0$ a subcomplex of
$BG_1{\cap}BG_2$. Choose $B{\Gamma}$ to be the double mapping cylinder of the inclusion maps.
Then there is a natural map ${\rho}: B{\Gamma} \to I$, where $I$ is the unit interval.
Let $E{\Gamma}$ be the universal cover of $B{\Gamma}$. Let 
$$E{\Gamma}{\times}_{\Gamma}{\E}{\G}_{\Gamma} 
\;\xrightarrow{q_{\Gamma}}\; I {\times} {\B}{\G}_{\Gamma}$$
be maps induced by the natural projection. 

We work as in Section \ref{sec-circle}.
Choose coset representatives:
$${\Delta}_i = \{{\gamma}_{i,j}: \;{\gamma}_{i,j}\in{\A}_i\}, \; i = 0, 1, 2.$$
In other words,
$${\Gamma} = \coprod_{j\in{\A}_i} G_i{\gamma}_{i,j},\; i = 0, 1, 2.$$
With the above notation:
\begin{itemize}
\item Every simplex $\sigma$ of ${\E}{\G}_{\Gamma}$ defines an equivalence relation on 
${\Delta}_i$,
$${\gamma}_{i,j} {\sim}_{\sigma}  {\gamma}'_{i,j} 
\;\Longleftrightarrow\; {\Gamma}_{\sigma}{\cap}{\gamma}_{i,j}G_i({\gamma}'_{i,j})^{-1} \not=
\emptyset.$$
\item Let $t\in I$. Set $i = 0$ if $t \in\text{Int}(I)$ and $i = t$ if $t\in{\partial}I$. For
each ${\gamma}_{i,j}\in{\Delta}_i$, a map
$${\iota}_{{\gamma}_{i,j}}: EG_i{\times}_{{\Gamma}_{\sigma}{\cap}G_i}g_{i,j}{\hat{\sigma}} \to
EG_i{\Gamma}{\times}_{{\Gamma}}{\Gamma}{\hat{\sigma}}, \;
[x, {\gamma}_{i,j}\hat{\sigma}] \mapsto [x, {\gamma}_{i,j}\hat{\sigma}]$$
here $EG_i{\Gamma} \subset E{\Gamma}$. Notice that the image of ${\iota}_{{\gamma}_{i,j}}$
 depends on the choice of $t\in I$.
\end{itemize}

The analogues of Lemmata \ref{lem-inverse} and \ref{lem-induced} hold in this case. Thus there
is a homeomorphism:
$${\chi}_i: EG_i{\times}_{G_i}G_i{\gamma}_{i,j}\hat{\sigma} \to 
E{\Gamma}{\times}_{\Gamma}{\E}{\G}_{\Gamma}$$
onto $\text{Im}({\iota}_{{\gamma}_{i,j}})$.

Equip $I$ with the structure of a simplicial complex 
with one $1$-simplex $I$,  and two $0$-simplices $0, 1$. The projection map induces a
functor:
$${\rho}: \text{cat}(I{\times}B{\G}_{\Gamma})^{op} \to \text{cat}(I)^{op}$$
We will use Segal's Pushdown Construction. We start by giving a description
of the functor $P_*{\rho}$ on $\text{cat}(I)^{op}$. For each simplex $t$ of $I$,
define a functor:
$${\bbS}{\circ}\text{bar}(q_{\Gamma})(t, -): \text{cat}(B{\G}_{\Gamma})^{op} \to 
{\bbS}pectra, \; {\sigma} \mapsto {\bbS}{\circ}\text{bar}(q_{\Gamma})(t, {\sigma})$$
Then the functor $P_*{\rho}$ is defined as:
$$P_*{\rho}: \text{cat}(I)^{op} \to {\bbS}pectra, \;
P_*{\rho}(t) = \text{hocolim}_{\text{cat}(B{\G}_{\Gamma})^{op}}
({\bbS}{\circ}\text{bar}(q_{\Gamma})).$$
So if we fix $t$, the functor $\text{bar}(q_{\Gamma})$ associates to $\sigma$, the space
$$q_{\Gamma}^{-1}(\hat{t}{\times}\hat{\sigma}) = (EG_i){\Gamma}{\times}_{\Gamma}
{\Gamma}\hat{\sigma} 
\cong (EG_i){\Gamma}{\times}_{{\Gamma}_{\sigma}}\hat{\sigma},$$ after subdividing,
where $i = 0$ if $t$ is an interior point and $G_i = G_{t+1}$ when $t\in{\partial}I$.
Actually, the inverse image for $t\in\text{Int}(I)$, is
$$q_{\Gamma}^{-1}(\hat{t}{\times}\hat{\sigma}) = (EG_i{\times}\{t\}){\Gamma}{\times}_{\Gamma}
{\Gamma}\hat{\sigma}$$
The inclusion map $G_i \to \Gamma$ induces a commutative diagram:
$$\begin{CD}
EG_i{\times}_{G_i}{\E}{\G}_{\Gamma} @>{f}>> 
EG_i{\Gamma}{\times}_{\Gamma}{\E}{\G}_{\Gamma} \\
@V{p_G}VV @VV{q_{\Gamma}}V \\
{\E}{\G}_{\Gamma}/G_i = {\B}{\G}_{G_i} @>{u}>> {\hat t}{\times}{\B}{\G}_{\Gamma} = 
{\hat t}{\times}{\E}{\G}_{\Gamma}/{\Gamma}
\end{CD}$$
where $u$ is just the quotient map. The analogue of Lemma \ref{lem-basic} also works here.

\begin{lem}\label{lem-basic2}
Suppose that the $\bbS$-IC holds for $G_i$, $i = 0, 1, 2$.
For each $t \in \text{cat}(I)^{op}$, there is a natural homotopy equivalence:
$$P_*{\rho}(t) = \text{hocolim}_{\text{cat}(B{\G}_{\Gamma})^{op}}
({\bbS}{\circ}\text{bar}(q_{\Gamma})(t, -)) \simeq {\bbS}(BG_i).$$
\end{lem}

The analogue of Proposition \ref{prop-circle} works in this case too.

\begin{prop}\label{prop-intrerval}
Assume that the ${\bbS}$-IC holds for $G_i$, $i = 0, 1, 2$. Then the following is a homotopy
cartesian diagram:
$$\begin{CD}
{\bbS}(BG_0) @>>> {\bbS}(BG_1) \\
@VVV @VVV \\
{\bbS}(BG_1) @>>> {\bbH}.(I{\times}{\B}{\G}_{\Gamma}, {\bbS}(q_{\Gamma}))
\end{CD}$$
\end{prop}

\begin{proof}
The proof works as in Proposition \ref{prop-circle}. The only difference is that for the
final result \cite{mupr:wa}, Section 2 is used.
\end{proof}

We also have the analogue of Propositions \ref{prop-iso} and \ref{prop-epi} when $\bbS = 
{\bbK}_R$.

\begin{prop}\label{prop-epi-interval}
Let $R$ be a ring such that $RG_0$ is regular coherent. Then 
\begin{enumerate}
\item the forgetful map 
$$f: {\bbH}.(I{\times}{\B}{\G}_{\Gamma}: {\bbK}_R(q_{\Gamma})) \to {\bbK}_R(B{\Gamma})$$
is a homotopy equivalence.
\item the assembly map
$$A: {\bbH}.({\B}{\G}_{\Gamma}, {\bbK}_R(p_{\Gamma})) \to {\bbK}_R(B{\Gamma})$$
induces an epimorphism on homotopy.
\end{enumerate}
\end{prop}

\begin{proof}
The proof follows from the splitting theorem in \cite{wa:wh} and \cite{wa:al}. The assumption
on the ring guarantees that the Waldhausen's Nil- groups vanish.
\end{proof}

\section{Applications}

We will apply the results to prove that ${\bbK}_R$-IC is true for certain types of groups. Before
we start we state the Novikov Conjecture for a discrete group $\Gamma$. Let ${\C}_{\Gamma}$
be the classifying space for the class of finite subgroups of $\Gamma$. Again, there is
a commutative diagram:
$$\begin{CD}
E{\Gamma}{\times}_{\Gamma}{\E}{\C}_{\Gamma} @>>> B{\Gamma} \\
@V{p}VV @VVV \\
{\B}{\C}_{\Gamma} @>>> *
\end{CD}$$
A group $\Gamma$ satisfies the {\bf integral $\bbS$-Novikov Conjecture} if the assembly map
$$A_{\C}: {\bbH}.({\B}{\C}_{\Gamma}, {\bbS}(p)) \to {\bbS}(B{\Gamma})$$
induces a split injection on the homotopy groups. 
It is an open question if all the torsion free groups 
satisfy the integral ${\bbK}_R$-Novikov Conjecture. It was proved that groups of finite
cohomological dimension that also have  finite asymptotic dimension satisfy the integral
${\bbK}_R$-Novikov Conjecture (\cite{bar}), generalizing the calculations in \cite{cape}.

A connection between the integral Novikov conjecture and the Isomorphism conjecture is
given in the following statement.

\begin{lem}\label{lem-nov}
Let $R$ be a regular coherent ring. Let $G$ be a torsion free group that satisfies the 
integral ${\bbK}_R$-Novikov Conjecture. Then the assembly map
$$A_{\C}: {\bbH}.({\B}{\C}_G, {\bbK}_R(p)) \to {\bbK}_R(B{\Gamma})$$
induces a monomorphism on the homotopy groups.
\end{lem}

\begin{proof}
Following the ideas in the Appendix in \cite{fajo:is}, notice that 
$A_{\C} = A{\circ}A_{{\C},{\G}}$ because of the commutative diagram
$$\begin{CD}
EG{\times}_G{\E}{\C}_G @>>> EG{\times}_G{\E}{\G}_G @>>> BG \\
@V{p}VV @VV{p_G}V @VVV \\
{\B}{\C}_G @>>> {\B}{\G}_G @>>> *
\end{CD}$$
The relative assembly map is the map induced by the first commutative square. Theorem A.10 in
\cite{fajo:is} states that $A_{{\C},{\G}}$ is an equivalence if the assembly map $A_{{\C}_S}$
is an equivalence for all virtually cyclic subgroups $S$ of $G$ (a similar calculation appears
in \cite{copr}). 
Since $G$ is torsion free, the
only finite subgroup of $G$ is the trivial group and the only virtually cyclic subgroups are
infinite cyclic subgroups. Thus $A_{{\C}_S}$ is the assembly map in the integral ${\bbK}_R$-Novikov
conjecture for $S \cong {\bbZ}$. Since $R$ is regular coherent, this map is an equivalence
(\cite{fajo:is}, Remark A.11). 
Therefore, since
$A_{\C}$ induces a monomorphism on homotopy groups, so does $A$.
\end{proof}

\begin{thm}\label{thm-iso}
Let $G$ be a torsion free group and $R$ a regular coherent ring. Assume that 
$G$ satisfies the ${\bbK}_R$-IC, and that $RG$ is a regular ring. 
Let $\Gamma$ be a torsion free group defined by:
\begin{itemize}
\item[(i)] ${\Gamma} = G{\rtimes}{\bbZ}$, or
\item[(ii)] ${\Gamma}= G_1*_GG_2$, such that $G_i$, $i = 1, 2$, satisfy the ${\bbK}_R$-IC,
\end{itemize}
such that $\Gamma$ satisfies the integral ${\bbK}_R$-Novikov conjecture. 
Then $\Gamma$ satisfies the ${\bbK}_R$-IC.
\end{thm}

\begin{proof}
Using Lemma \ref{lem-nov}, we see that the assembly map $A$ induces a monomorphism on the
homotopy groups. Using Proposition \ref{prop-epi} for case (i) and Proposition 
\ref{prop-epi-interval} for case (ii), we see that $A$ is an epimorphism. Thus $A$ is an
equivalence.
\end{proof}

Now we give applications of the Theorem.

\begin{cor}\label{cor-asym}
Let $G$ be a torsion free group and $R$ a regular coherent ring. Assume that:
\begin{enumerate} 
\item $G$ satisfies the ${\bbK}_R$-IC,
\item $G$ has finite asymptotic dimension,
\item $RG$ is a regular ring. 
\end{enumerate}
Let $\Gamma$ be a torsion free group defined by:
\begin{itemize}
\item[(i)] ${\Gamma} = G{\rtimes}{\bbZ}$, or
\item[(ii)] ${\Gamma}= G_1*_GG_2$, such that $G_i$, $i = 1, 2$, satisfy the ${\bbK}_R$-IC and
they have finite cohomological and asymptotic dimensions.
\end{itemize}
Then $\Gamma$ satisfies the ${\bbK}_R$-IC.
\end{cor}

\begin{proof}
The assumptions on the groups imply that $\Gamma$ has finite asymptotic dimension (\cite{bedr}). 
Also, $\Gamma$ has finite cohomological dimension and thus it admits a finite dimensional 
$B{\Gamma}$. By \cite{bar}, $\Gamma$ satisfies the integral ${\bbK}_R$-Novikov conjecture.
The result follows from Theorem \ref{thm-iso}.
\end{proof}

\begin{cor}\label{cor-free}
Let $F$ be a finitely generated free group and $R$ a regular coherent ring. Then 
\begin{enumerate}
\item $F$ satisfies the ${\bbK}_R$-IC.
\item $F{\rtimes}\bbZ$ satisfies the ${\bbK}_R$-IC
\end{enumerate}
\end{cor}

\begin{proof}
The first statement follows by induction on the number $k$ of generators of $F$:
If $k = 2$, then $F = {\bbZ}*{\bbZ}$ and the result follows from \ref{cor-asym}. 
For $k > 2$, $F = F_{k-1}*{\bbZ}$, where $F_{k-1}$ is the free group on $(k - 1)$ 
generators. All the assumptions of Corollary \ref{cor-asym} are satisfied and thus
$F$ satisfies the ${\bbK}_R$-IC. 

For the second statement, we use again Corollary \ref{cor-asym}. Part (1) implies that
assumption (1) is satisfied. Also, the free group has finite asymptotic dimension and thus
assumption (2) is satisfied. Assumption (3) follows because $R$ is regular Noetherian ring.
\end{proof}

\begin{rem}
In \cite{arfarou} and \cite{farou}, there was a special assumption for groups of type
$F{\rtimes}{\bbZ}$ to satisfy the fibered pseudoisotopy IC. 
Essentially the assumption was that the action of $\bbZ$ on $F$ has certain
geometric properties. 
Corollary \ref{cor-free} is more general because such an assumption is not needed but it only
gives the ${\bbK}_R$-IC and not the fibered version. A general fibered version could not
follow along the same lines because the assumption in Corollary \ref{cor-asym} guarantee
that all the Nil-groups that appear vanish.
\end{rem}

\end{document}